\documentclass[12pt]{amsart}
\usepackage{amsmath,amsthm}
\usepackage{amssymb}
\usepackage{euscript}
\DeclareMathAlphabet{\EuRm}{U}{eur}{m}{n}
\SetMathAlphabet{\EuRm}{bold}{U}{eur}{b}{n}
\begin{document}
%
%
\swapnumbers
\newtheorem{thm}{Theorem}[section]
\newtheorem*{TA}{Theorem A}
\newtheorem*{TB}{Theorem B}
\newtheorem{lemma}[thm]{Lemma}
\newtheorem{prop}[thm]{Proposition}
\newtheorem{cor}[thm]{Corollary}
\theoremstyle{definition}
\newtheorem{assum}[thm]{Assumption}
\newtheorem{defn}[thm]{Definition}
\newtheorem{example}[thm]{Example}
\newtheorem{notation}[thm]{Notation}
\newtheorem{summary}[thm]{Summary}
\newtheorem{fact}[thm]{Fact}
\theoremstyle{remark}
\newtheorem{remark}[thm]{Remark}
\newtheorem{assume}[thm]{Assumption}
\newtheorem{note}[thm]{Note}
\newtheorem{ack}[thm]{Acknowledgements}

\numberwithin{equation}{section}
%
%
\def\sect{\setcounter{thm}{0} \section}
%
%
\newcommand{\xra}[1]{\xrightarrow{#1}}\
\newcommand{\xla}[1]{\xleftarrow{#1}}
\newcommand{\hra}{\hookrightarrow}
\newcommand{\adj}[2]{\substack{{#1}\\ \rightleftharpoons \\ {#2}}}
\newcommand{\hsp}{\hspace{10 mm}}
\newcommand{\hs}{\hspace{5 mm}}
\newcommand{\hsm}{\hspace{2 mm}}
\newcommand{\vs}{\vspace{7 mm}}
\newcommand{\vsm}{\vspace{2 mm}}
\newcommand{\rest}[1]{\lvert_{#1}}
\newcommand{\lra}[1]{\langle{#1}\rangle}
\newcommand{\DEF}{:=}
\newcommand{\EQUIV}{\Leftrightarrow}
\newcommand{\epic}{\to\hspace{-5 mm}\to}
\newcommand{\hotimes}{\hat{\otimes}}
\newcommand{\msim}{/\!\sim}
%
%
\newcommand{\Coker}{\operatorname{Coker}}
\newcommand{\colim}{\operatorname{colim}}
\newcommand{\diag}{\operatorname{diag}}
\newcommand{\Fib}{\operatorname{Fib}}
\newcommand{\Hom}{\operatorname{Hom}}
\newcommand{\Image}{\operatorname{Im}}
\newcommand{\Ker}{\operatorname{Ker}}
\newcommand{\sk}[1]{\operatorname{sk}_{#1}}
%
%
\newcommand{\Alg}{{\EuScript Alg}}
\newcommand{\C}{{\mathcal C}}
\newcommand{\sC}[2]{s_{\lra{#1}}{#2}}
\newcommand{\D}{{\mathcal D}}
\newcommand{\gD}[1]{{\EuScript D}^{#1}}
\newcommand{\F}{{\mathcal F}}
\newcommand{\G}{{\mathcal G}}
\newcommand{\Gp}{{\EuScript Gp}}
\newcommand{\I}{{\EuScript I}}
\newcommand{\Pa}{$\Pi$-algebra}
\newcommand{\PAlg}{\hbox{$\Pi$-$\Alg$}}
\newcommand{\Ss}{{\mathcal S}}
\newcommand{\Sa}{\Ss_{\ast}}
\newcommand{\gS}[1]{{\EuScript S}^{#1}}
\newcommand{\TT}{{\mathcal T}}
\newcommand{\Ta}{\TT_{\ast}}
\newcommand{\Tc}{\TT_{c}}
%
%
\newcommand{\DD}{{\mathbb \Delta}}
\newcommand{\N}{\mathbb N}
\newcommand{\bP}{\mathbb P}
\newcommand{\R}{\mathbb R}
\newcommand{\bk}{\mathbf{k}}
\newcommand{\bn}{\mathbf{n}}
%
%
\newcommand{\A}{\mathbf{A}}
\newcommand{\B}{\mathbf{B}}
\newcommand{\De}[1]{{\mathbf \Delta}[{#1}]}
\newcommand{\dD}[1]{\overset{\bullet}{\mathbf{\Delta}}[{#1}]}
\newcommand{\bD}[1]{\mathbf{D}^{#1}}
\newcommand{\hD}[1]{\hat{\mathbf{D}}^{#1}}
\newcommand{\bS}[1]{\mathbf{S}^{#1}}
\newcommand{\hS}[1]{\hat{\mathbf{S}}^{#1}}
\newcommand{\Pe}[1]{\mathbf{P}^{#1}}
\newcommand{\PP}[2]{\Pe{#1}_{#2}}
\newcommand{\U}{\mathbf U}
\newcommand{\V}{\mathbf V}
\newcommand{\W}{\mathbf W}
\newcommand{\X}{\mathbf{X}}
\newcommand{\Y}{\mathbf{Y}}
\newcommand{\Z}{\mathbf{Z}}
%
%
\newcommand{\pis}{\pi_{\ast}}
\newcommand{\Gs}{G_{\ast}}
\newcommand{\Hs}{H_{\ast}}
\newcommand{\Js}{J_{\ast}}
\newcommand{\Ks}{K_{\ast}}
\newcommand{\Ls}{L_{\ast}}
\newcommand{\Ts}{T_{\ast}}
%
%
\newcommand{\Ad}{A_{\bullet}}
\newcommand{\Bd}{B_{\bullet}}
\newcommand{\Edd}{(E_{\bullet})^{\bullet}_{\Delta}}
\newcommand{\Eud}{E_{\bullet}^{\bullet}}
\newcommand{\Fd}{F^{\bullet}}
\newcommand{\Fdd}{\Fd_{\Delta}}
\newcommand{\Kd}{K_{\bullet}}
\newcommand{\Ld}{L_{\bullet}}
\newcommand{\Qd}{Q_{\bullet}}
\newcommand{\Rd}{R_{\bullet}}
\newcommand{\Ud}{\U^{\bullet}}
\newcommand{\Udd}{\U^{\bullet}_{\Delta}}
\newcommand{\Vd}{\V_{\bullet}}
\newcommand{\Vdd}{(\V_{\bullet})^{\bullet}_{\Delta}}
\newcommand{\Wd}{\W_{\bullet}}
\newcommand{\uW}[2]{\W^{#1}_{#2}}
\newcommand{\bW}[2]{\bar{\W}^{#1}_{#2}}
\newcommand{\uWd}[1]{\uW{#1}{\bullet}}
\newcommand{\Wdd}{(\Wd)^{\bullet}_{\Delta}}
\newcommand{\Wud}{\uWd{\bullet}}
\newcommand{\Xd}{X_{\bullet}}
\newcommand{\Zud}{\Zd_{\bullet}^{\bullet}}

\newcommand{\co}[1]{c{#1}_{\bullet}}
%
%
%
\title[CW resolutions of spaces]{CW simplicial resolutions of spaces\\
with an application to loop spaces}
\author{David Blanc}
\address{Dept.\ of Mathematics, University of Haifa, 31905 Haifa, Israel}
\email{blanc@math.haifa.ac.il}
\date{April 23, 1998}
\begin{abstract}
We show how a certain type of $CW$ simplicial resolutions of spaces by wedges
of spheres may be constructed, and how such resolutions yield an obstruction 
theory for a given space to be a loop space.
\end{abstract}

\maketitle
%
%
\sect{Introduction}
\label{ci}

A simplicial resolution of a space $\X$ by wedges of spheres is a simplicial 
space \ $\Wd$ \ such that \ (a) each space \ $\W_{n}$ \ is
homotopy equivalent to a wedge of spheres, and \ (b) for each \ $k\geq 1$, \ 
the augmented simplicial group \ $\pi_{k}\Wd\to\pi_{k}\X$ \ is acyclic 
(see \S \ref{dfsr} below). Such resolutions, first constructed by Chris Stover 
in \cite[\S 2]{StoV}, are dual to the ``unstable Adams resolutions'' of 
\cite[I, \S 2]{BKaH}, and have a number of applications:
see \S \ref{sqss} below and \cite{StoV,DKSStP,DKStE,BlaH,BlaHH,BlaHO,BlaL}. 

However, the Stover construction yields very large resolutions, which do not 
lend themselves readily to computation, and no other construction was hitherto
available.  In particular, it was not clear whether one could find 
minimal resolutions of this type. The purpose of this note is to show that 
any space $X$ has simplicial resolutions by wedges of spheres, which may be 
constructed from purely algebraic data, consisting of an (arbitrary) 
simplicial resolution of \ $\pis X$ \ as a \Pa \ -- \ that is, as
a graded group with an action on the primary homotopy operations on it 
(see \S \ref{dpa} below):

\begin{TA}
Every free simplicial \Pa\ resolution of a realizable \Pa\ \ $\pis X$ \ 
is realizable topologically as a simplicial resolution by wedges of spheres.
\end{TA}

\noindent and in fact such resolutions can be given a convenient ``CW structure'' 
(\S \ref{dcwg}). \ There is an analogous result for maps (Theorem \ref{ttwo}). 

Since \emph{no} such resolution of a 
non-realizable \Pa\ can be realized (see \S \ref{rpa} below), this completely 
determines which free simplicial \Pa\ resolutions are realizable.

The Theorem implies that in the spectral sequences of \cite{StoV,BlaH,DKSStP} 
we can work with minimal resolutions, and allows us to identify the 
higher homotopy operations of \cite{BlaHH,BlaH,BlaL} as lying in appropriate
cohomolgy groups (compare \cite[4.17]{BlaHO} and \cite[\S 6]{BlaHR}). A
generalization of Theorem A to other model categories appears in 
\cite{BlaAI}\vsm .

As an application of such CW resolutions, we describe an obstruction theory 
for deciding whether a given space $\X$ is a loop space, 
in terms of higher homotopy operations. One such theory was given in 
\cite{BlaL}, but the present approach does not require a given $H$-space 
structure on $\X$, and may be adapted also to the existence of \ 
$A_{n}$-structures (and thus subsumes \cite{BlaHO}). It is summarized in

\begin{TB}
A space $\X$ with trivial Whitehead products is homotopy equivalent to a 
loop space if and only if the higher homotopy operations of \S \ref{dhho} 
below vanish coherently.
\end{TB}

\subsection{Notation and conventions}
\label{snac}\stepcounter{thm}

$\Gp$ \ will denote the category of groups, \ $\TT$ \ that of topological 
spaces, and $\Ta$ \ that of pointed topological spaces with base-point 
preserving maps. The full subcategory of $0$-connected spaces will be denoted 
by \ $\Tc\subset\Ta$. \ The category of simplical sets will be denoted by $\Ss$ \ 
and that of pointed simplicial sets by \ $\Sa$; \  we shall use boldface
letters: \ $\X$, \ $\bS{n},\dotsc$ \ to denote objects in any of these four 
categories. If \ $f:\X\to\Y$ \ is a map in one of these categories, we denote
by \ $f_{\#}:\pis\X\to\pis\Y$ \ the induced map in the homotopy groups.

\subsection{Organization}
\label{sorg}\stepcounter{thm}

In section \ref{cso} we review some background on simplicial objects and 
bisimplicial groups, and in section \ref{cpar} we recall some facts on \Pa s, 
and prove our main results on CW resolutions of spaces by wedges of spheres: 
Theorem A (=Theorem \ref{tone}) and Theorem \ref{ttwo}.
In section \ref{csbc} we define a certain
cosimplicial simplicial space up-to-homotopy, which can be rectified
if and only if $\X$ is a loop space. In section \ref{cphho} we construct a
certain collection of \textit{face-codegeneracy polyhedra}, which are used
to define the higher homotopy operations refered to in Theorem B
(=Theorem \ref{tthree}). We also show how the theorem may be used in reverse 
to calculate a certain tertiary operation in \ $\pis\bS{7}$.

\begin{ack}\label{aa}\stepcounter{subsection}
I would like to thank the referee for his or her comments (see in particular 
\S \ref{rref} below).
\end{ack}
%
%
%
\sect{simplicial objects}
\label{cso}

We first provide some definitions and facts on simplicial objects:

\begin{defn}\label{dso}\stepcounter{subsection}
Let \ $\DD$ \ denote the category of ordered sequences \ 
$\bn= \lra{0,1,\dotsc,n}$ \ ($n\in \N$), \ with 
order-preserving maps. 
A \textit{simplicial object} over a category $\C$  is a functor \ 
$X:\DD^{op}\to\C$, \ usually written \ $\Xd$, \ which may be described 
explicitly as a sequence of objects \ $\{ X_{k} \}_{k=0}^{\infty}$ \ in $\C$, 
equipped with face maps \ $d^{k}_{i}:X_{k}\to X_{k-1}$ \ and degeneracies \ 
$s^{k}_{j}:X_{k}\to X_{k+1}$ \ (usually written simply \ $d_{i}$, \ $s_{j}$, \ 
for \ $0\leq i,j\leq k$), \ satisfying the usual 
simplicial identities (\cite[\S 1.1]{MayS}). \ 
If \ $I=(i_{1},i_{2},\dotsc,i_{r})$ \ is some multi-index, \ we write \ 
$d_{I}$ \ for \ $d_{i_{1}}\circ d_{i_{2}}\circ\dotsb\circ d_{i_{r}}$, \ with \ 
$d_{\emptyset}\DEF id$; \ and similarly for \ $s_{I}$. \ An \textit{augmented}
simplicial object is one equipped with an augmentation \ 
$\varepsilon:X_{0}\to Y$ \ (for \ $Y\in\C$), \ with \ 
$\varepsilon d_{0}=\varepsilon d_{1}$.

The category of simplicial objects over $\C$ is denoted by \ $s\C$.
We write \ $\sC{n}{\C}$ \ for the category $n$-\textit{simplicial objects} 
over $\C$ (that is, objects of the form \ $\{ X_{k} \}_{k=0}^{n}$, \ with
the relevant face maps and degeneracies), \ and denote the
truncation functor \ $s\C\to\sC{n}{\C}$ \ by \ $\tau_{n}$.
\end{defn}

For technical convenience in the next two sections we shall be working 
mainly in the category of simplicial groups, denoted by \ $\G$ \ (rather 
than \ $s\Gp$); \ objects in $\G$ will be denoted by capital letters $X$, $Y$, 
and so on. \ A simplicial object \ $\Xd=(X_{0},X_{1},\dotsc)$ \ in \ 
$s\G$ \ is thus a bisimplicial group, which has an \textit{external} 
simplicial dimension (the $n$ in \ $X_{n}\in\G$), \ as well as the 
\textit{internal} simplicial dimension $k$ (inside $\G$), \ which 
we shall denote by \ $(X_{n})_{k}^{int}$, \ if necessary.

\subsection{Simplicial sets and groups}
\label{ssg}\stepcounter{thm}

The standard $n$ simplex in $\Ss$ is denoted by \ 
$\De{n}$, \ generated by \ $\sigma_{n}\in \De{n}_{n}$. \ 
$\dD{n}$ \ denotes the sub-object of \ $\De{n}$ \ generated by \ 
$d_{i}\sigma_{n}$ \ ($0\leq i\leq n$). \ The simplicial $n$-sphere is \ 
$\bS{n}\DEF \De{n}/\dD{n}$, \ and the $n$-disk is \ $\bD{n}\DEF C\bS{n-1}$.

Let \ $F:\Ss\to\G$ \ denote the (dimensionwise) free group functor of
\cite[\S 2]{MilnF}, and \ $G:\Ss\to\G$ \ be Kan's simplicial loop functor \ 
(cf.\ \cite[Def.\ 26.3]{MayS}), \ with \ $\bar{W}:\G\to\Ss$ \ the 
Eilenberg-Mac Lane classifying space functor (cf.\ \cite[\S 21]{MayS}). \ 
Recall that if \ $S:\TT\to\Ss$ \ is the singular set functor and \ 
$\|-\|:\Ss\to\TT$ \ the geometric realization functor 
(see \cite[\S 1,14]{MayS}), then the adjoint pairs of functors
%
\setcounter{equation}{\value{thm}}\stepcounter{subsection}
\begin{equation}\label{ezero}
\TT\ \ \substack{S\\ \rightleftharpoons\\ \|-\|}\ \ \Ss\ \ 
\substack{G\\ \rightleftharpoons\\ \bar{W}}\ \ \G
\end{equation}
\setcounter{thm}{\value{equation}}
induce isomorphisms of the corresponding homotopy categories (see
\cite[I, \S 5]{QuH}), so that for the purposes of homotopy theory we can work in 
$\G$ rather than $\TT$.

\begin{defn}\label{dgs}\stepcounter{subsection}
In particular, \ $\gS{n}\DEF F\bS{n-1}\in\G$ \ for \ $n\geq 1$ \ (and \ 
$\gS{0}\DEF G\bS{0}$ \ for \ $n=0$) \ 
will be called the \textit{$n$-dimensional $\G$-sphere},  in as much as \ 
$[\gS{n},G\X]_{\G}\cong\pi_{n}\X=[\bS{n},\X]$ \ for any Kan 
complex \ $\X\in\Ss$. \ Similarly, \ 
$\gD{n}\DEF F\bD{n-1}$ \ will be called the \textit{$n$-dimensional $\G$-disk}.
\end{defn}

\begin{defn}\label{dmat}\stepcounter{subsection}
In any complete category $\C$, the \textit{matching object} 
functor \ $M:\Ss^{op}\times s\C\to\C$, \ 
written \ $M_{\A}\Xd$ \ for \ a (finite) simplicial set \ $\A\in\Ss$ \ and \ 
$\Xd\in s\C$, \ is defined by requiring: \ (a) \ $M_{\De{n}}\Xd:=X_{n}$, \ 
and \ (b) \ if \ $\A=\colim_{i}\A_{i}$, \ then \ 
$M_{\A}\Xd=\lim_{i} M_{\A_{i}}\Xd$ \ (see \cite[\S 2.1]{DKStB}). \ In 
particular, if \ $\A^{k}_{n}$ \ is the subcomplex of \ $\dD{n}$ \ generated by 
the last \ $(n-k+1)$ \ faces \ $(d_{k}\sigma_{n},\dotsc,d_{n}\sigma_{n})$, \ 
we write \ $M^{k}_{n}\Xd$ \ for \ $M_{\A^{k}_{n}}\Xd$: \ explicitly,
%
\setcounter{equation}{\value{thm}}\stepcounter{subsection}
\begin{equation}\label{eone}
M^{k}_{n}\Xd=\{(x_{k},\dotsc,x_{n})\in(X_{n-1})^{n+1}~\lvert\ 
d_{i}x_{j}=d_{j-1}x_{i} \text{\ \ for all\ }k\leq i<j\leq n\}.
\end{equation}
\setcounter{thm}{\value{equation}}

\noindent and the map \ $\delta^{k}_{n}:X_{n}\to M^{k}_{n}\Xd$ \ induced by 
the inclusion \ $\A^{k}_{n}\hra\De{n}$  \ is defined \ 
$\delta^{k}_{n}(x)=(d_{k}x,\dotsc,d_{n}x)$. \ The original matching object of
\cite[X,\S 4.5]{BKaH} \ was \ $M_{n}^{0}\Xd=M_{\dD{n}}\Xd$, \ 
which we shall further abbreviate to \ $M_{n}\Xd$; \ each face map \ 
$d_{k}:X_{n+1}\to X_{n}$ \ factors through \ $\delta_{n}\DEF\delta^{0}_{n}$. \ 
See also \cite[XVII, 87.17]{PHirL}.
\end{defn}

\begin{remark}\label{rmat}\stepcounter{subsection}\
Note that for \ $X\in \G$ \ and \ $\A\in\Ss$ \ we have \ 
$M_{\A}X\cong \Hom_{\G}(F\A,X)\in\Gp$ \ (cf.\ \S \ref{ssg}), so for \ 
$\Xd\in s\G$ \ also \ $(M_{\A}X)_{k}\cong \Hom_{\G}(F\A,(\Xd)^{int}_{k})$ \ in 
each simplicial dimension $k$.
\end{remark}

\begin{defn}\label{dfib}\stepcounter{subsection}
$\Xd\in s\G$ \ is called \textit{fibrant} if each of the maps \ 
$\delta_{n}:X_{n}\to M_{n}\Xd$ \ ($n\geq 0$) \ is a fibration in $\G$ \ 
(that is, a surjection onto the identity component \ -- \ see 
\cite[II, 3.8]{QuH}). This is just the condition for fibrancy in the Reedy
model category, (see \cite{ReedH}), as well as in that of \cite{DKStE}, but 
we shall not make explicit use of either.
\end{defn}

By analogy with Moore's normalized chains (cf.\ \cite[17.3]{MayS}) we have:

\begin{defn}\label{dnc}\stepcounter{subsection}
Given \ $\Xd\in s\G$, \ we define the 
$n$-\textit{cycles} object of \ $\Xd$, \ written \ $Z_{n}\Xd$, \ to be the 
fiber of \ $\delta_{n}:X_{n}\to M_{n}\Xd$, \ so \ 
$Z_{n}\Xd=\{ x\in X_{n}\,|\ d_{i}x=0 \ \text{for}\ i=0,\dotsc,n\}$ \ 
(cf.\ \cite[I,\S 2]{QuH}). \ 
Of course, this definition really makes sense only when \ $\Xd$ \ is fibrant 
(\S \ref{dfib}). \ Similarly, the $n$-\textit{chains} object
of \ $\Xd$, \ written \ $C_{n}\Xd$, \ is defined to be the fiber of \ 
$\delta^{1}_{n}:X_{n}\to M^{1}_{n}\Xd$.
\end{defn}

If \ $\Xd\in s\G$ \ is fibrant, the map \ 
$d_{0}^{n}=d_{0}\rest{C_{n}\Xd}:C_{n}\Xd\to Z_{n-1}\Xd$ \ is the pullback of \ 
$\delta_{n}:X_{n}\to M_{n}\Xd$ \ along the inclusion \ 
$\iota:Z_{n-1}\Xd\to M_{n}\Xd$ \ (where \ $\iota(z)=(z,0,\dotsc,0)$), \ so \ 
$d_{0}^{n}$ \ is a fibration (in $\G$), \ fitting into a fibration sequence
%
\setcounter{equation}{\value{thm}}\stepcounter{subsection}
\begin{equation}\label{etwo}
Z_{n}\Xd\xra{j_{n}}C_{n}\Xd \xra{d^{n}_{0}} Z_{n-1}\Xd.
\end{equation}
\setcounter{thm}{\value{equation}}

%
%
\begin{prop}\label{pone}\stepcounter{subsection}
For any fibrant \ $\Xd\in s\C$, \ the inclusion \ $\iota:C_{n}\Xd\hra X_{n}$ \ 
induces an isomorphism \ 
$\iota_{\star}:\pis C_{n}\Xd\cong C_{n}(\pis\Xd)$ \ for each \ $n\geq 0$.
\end{prop}

\begin{proof}
(a)\hs First note that if \ $j:\A\hra\B$ \ is a trivial cofibration in $\Ss$, \ 
then \ $j^{\ast}:M_{\B}\Xd\to M_{\A}\Xd$ \ has a natural section \ 
$r:M_{\A}\Xd\to M_{\B}\Xd$ \ (with \ $j^{\ast}\circ r=id$) \ for any \ 
$\Xd\in s\G$: \ \ This is because by remark \ref{rmat}, \ 
$(M_{\A}\Xd)_{k}\cong\Hom_{\G}(F\A,(\Xd)^{int}_{k})$ \ for \ $\A\in\Ss$; \ 
since \ $F\A$ \ is fibrant in $\G$, we can choose a left inverse \ 
$\rho:F\B\to F\A$ \ for  \ $Fj:F\A\hra F\B$, \ so \ 
$j^{\ast}:(M_{\B}\Xd)^{int}_{k}\to (M_{\A}\Xd)^{int}_{k}$ \ 
has a right inverse \ $\rho^{\ast}$, \ which is natural in \ 
$(\Xd)^{int}_{k}$; \ so these maps \ $\rho^{\ast}$ \ fit together to yield 
the required map $r$.

This need not be true in general if $j$ is not a weak 
equivalence, as the example of \ $M_{2}^{1}\Xd\to M_{1}^{0}\Xd$ \ shows.\vsm 

(b)\hs Given \ $\eta\in C_{n}\pi_{m}\Xd$ \ represented by \ 
$h:\gS{m}\to X_{n}$ \ with \ $d_{k}h\sim 0$ \ ($1\leq k\leq n$), \ consider 
the diagram:

%
%
\begin{picture}(250,180)(-60,-15)
\put(20,138){$\gS{m}$}
\put(32,135){\vector(1,-1){32}}
\put(57,114){$h$}
\put(35,140){\line(4,-1){112}}
\put(147,112){\vector(1,-1){40}}
\put(90,130){$\sim 0$}
\put(23,133){\line(1,-4){20}}
\put(43,53){\vector(4,-3){60}}
\put(18,60){$\sim 0$}
\put(65,95){$X_{n}$}
\put(80,90){\vector(1,-1){20}}
\put(93,80){$\delta^{k}_{n}$}
\put(84,100){\vector(3,-1){95}}
\put(110,95){$\delta^{k+1}_{n}$}
\put(70,90){\vector(1,-2){39}}
\put(75,47){$d_{k}$}
%
%
\put(100,60){$M^{k}_{n}\Xd$}
\put(135,65){\vector(1,0){43}}
\put(181,60){$M^{k+1}_{n}\Xd$}
\put(120,55){\vector(0,-1){42}}
\put(106,37){$\pi_{k}$}
\put(105,0){$X_{n-1}$}
\put(195,55){\vector(0,-1){42}}
\put(198,32){$j^{\ast}=(d_{k},\dotsc,d_{k})$}
\put(135,5){\vector(1,0){42}}
\put(148,12){$\delta^{k}_{n-1}$}
\put(181,0){$M^{k}_{n-1}\Xd$}
\put(133,43){$\framebox{\scriptsize{PB}}$}
\end{picture}

\noindent in which \ $j^{\ast}$ \ is a fibration by (a) if \ $k\geq 1$, \ 
so the lower left-hand square is in fact a homotopy pullback square 
(see \cite[\S 1]{MatP}). 
By descending induction on \ $1\leq k\leq n-1$, \ (starting with \ 
$\delta^{n}_{n}=d_{n}$), \ we may assume \ 
$\delta^{k+1}_{n}\circ h:\gS{m}\to M^{k+1}_{n}\Xd$ \ is nullhomotopic in $\C$,
as is \ $d_{k}\circ h$, \ so the induced pullback map \ 
$\delta^{k}_{n}\circ h:\gS{m}\to M^{k}_{n}\Xd$, \ is also nullhomotopic by the
universal property. We conclude that \ $\delta^{1}_{n}\circ h\sim 0$, \ 
and since \ $\delta^{1}_{n}:X_{n}\to M^{1}_{n}\Xd$ \ is a fibration by (a), 
we can choose \ $h:\gS{m}\to X_{n}$ \ so that \ $\delta^{1}_{n} h=0$. \ 
Thus \ $h$ \ lifts to \ $C_{n}\Xd=\Fib(\delta^{1}_{n})$, \ and \  
$\iota_{\star}$ is surjective\vsm .

(c)\hs Finally, the long exact sequence in homotopy for the fibration 
sequence \ 
$$
C_{n}\Xd\xra{\iota} X_{n}\xra{\delta^{1}_{n}} M^{1}_{n}\Xd
$$
implies that \ $\iota_{\#}:\pis C_{n}\Xd\to\pis X_{n}$ \ is monic, \ so \ 
$\iota_{\star}:\pis C_{n}\Xd\to C_{n}(\pis\Xd)$ \ is, too.
\end{proof}

\begin{defn}\label{dlat}\stepcounter{subsection}
The dual construction to that of \S \ref{dmat} yields the colimit \ 
$$
L_{n}\Xd\DEF\coprod_{0\leq i\leq n-1} X_{n-1}\msim,
$$

\noindent where for any \ 
$x\in X_{n-2}$ \ and \ $0\leq i\leq j \leq n-1$ \ we set \ 
$s_{j}x$ \ in the $i$-th copy of \ $X_{n-1}$ \ equivalent under $\sim$ to 
$s_{i}x$ \ in the $(j+1)$-st copy of \ $X_{n-1}$. \ $L_{n}\Xd$ \ 
has sometimes been called the ``$n$-th latching object'' of \ 
$\Xd$. \ The map \ $\sigma_{n}:L_{n}\Xd\to X_{n}$ \ is defined \ 
$\sigma_{n}x_{(i)}=s_{i}x$, \ where \ $x_{(i)}$ \ is in the $i$-th copy of \ 
$X_{n-1}$.
\end{defn}

%
%
\sect{\Pa s and resolutions}
\label{cpar}

In this section we recall some definitions and prove our main results on 
\Pa s and resolutions:

\begin{defn}\label{dpa}\stepcounter{subsection}
A \textit{\Pa} \ is a graded group \ $\Gs=\{G_{k}\}_{k=1}^{\infty }$ \ 
(abelian in degrees $>1$), together with an action on \ $\Gs$ \ of the 
primary homotopy operations (i.e., compositions and Whitehead products, 
including the ``$\pi_{1}$-action'' of \ $G_{1}$ \ on the higher \ 
$G_{n}$'s, \ as in \cite[X, \S 7]{GWhE}), \ satisfying the usual universal 
identities. See \cite[\S 2.1]{BlaA} for a more explicit description. \ These 
are algebraic models of the homotopy groups \ $\pis\X$ \ of a space 
(or Kan complex) $\X$, in the same way that an algebra over the Steenrod 
algebra models its cohomology ring. \ 
The category of \Pa s is denoted by \ $\PAlg$.

We say that a space (or Kan complex, or simplicial group) \ $\X$ 
\textit{realizes} an (abstract) \Pa \ $\Gs$ \ if there is an isomorphism of 
\Pa s \ $\Gs\cong \pis\X$. \ (There may be non-homotopy equivalent spaces 
realizing the same \Pa \ -- \ cf.\ \cite[\S 7.18]{BlaHH}).\hsm Similarly, 
an abstract morphism of \Pa s \ $\phi:\pis \X\to\pis\Y$ \ (between realizable 
\Pa s) is \textit{realizable} if there is a map \ $f:\X\to\Y$ \ such that \ 
$\pis f=\phi$.
\end{defn}

\begin{defn}\label{dfpa}\stepcounter{subsection}
The \textit{free \Pa\ generated} by a graded set \ 
$\Ts=\{T_{k}\}_{k=1}^{\infty }$ \ is \ $\pis\W$, \ where \ 
$\W=\bigvee_{k=1}^{\infty }\bigvee_{{\tau\in T_{k}}} \bS{k}_{(\tau)}$ \ (and we
identify \  $\tau\in T_{k}$ \ with the generator of \
$\pi_{k}\W $ \ representing the inclusion\ \ $\bS{k}_{(\tau)}\hra\W$). \ 

If we let \ $\F\subset\PAlg$ \ denote the full subcategory of free \Pa s, and
$\Pi $ the homotopy category of wedges of spheres 
(inside \ $ho\Ta$ \ or \ $ho\Sa$ \ -- \ or equivalently, the homotopy
category of coproducts of $\G$-spheres in \ $ho\G$), \ then
the functor \ $\pis:\Pi\to\F$ \ is an equivalence of categories. 
Thus any \Pa\ morphism \ $\varphi:\Gs\to\Hs$ \ is realizable (uniquely, up to
homotopy), if \ $\Gs$ \ and \ $\Hs$ \ are free \Pa s \ (actually, only \ $\Gs$ \ 
need be free).
\end{defn}

\begin{defn}\label{dfpc}\stepcounter{subsection}
Let \ $T:\PAlg\to\PAlg$ \ be the ``free \Pa'' comonad 
(cf.\ \cite[VI, \S 1]{MacC}), \ defined \ $T\Gs=
\coprod_{k=1}^{\infty}\coprod_{g\in G_{k}}~\pis\bS{k}_{(g)}$. \ 
The counit \ $\varepsilon=\varepsilon_{\Gs}:T\Gs\epic\Gs$ \ is defined by \ 
$\iota^{k}_{(g)}\mapsto g$ \ (where \ $\iota^{k}_{(g)}$ \ is the canonical 
generator of \ $\pis\bS{k}_{(g)}$), \ and the comultiplication \ 
$\vartheta=\vartheta_{\Gs}:T\Gs\hra T^{2}\Gs$ \ is induced by the 
natural transformation \ $\bar{\vartheta}:id_{\F}\to T|_{\F}$ \ 
defined by \ $x_{k}\mapsto \iota^{k}_{(x_{k})}$.
\end{defn}

\begin{defn}\label{dapa}\stepcounter{subsection}
An \textit{abelian} \Pa\ is one for which all Whitehead products vanish.
\end{defn}

These are indeed the abelian objects of \ $\PAlg$ \ -- \ see \cite[\S 2]{BlaA}.
In particular, if $\X$ is an $H$-space, then \ $\pis\X$ \ is an abelian \Pa\ 
(cf.\ \cite[X, (7.8)]{GWhE}).

\begin{defn}\label{dfsr}\stepcounter{subsection}
A simplicial \Pa\ \ $\Ad$ \ is called \textit{free} if for each \
$n\geq 0$ \ there is a graded set \ $\Ts^{n}\subseteq A_{n}$ \
such that \ $A_{n}$ \ is the free \Pa\ generated by \ $\Ts^{n}$ \ 
(\S \ref{dfpa}), and each degeneracy map \ $s_{j}:A_{n}\to A_{n+1}$ \
takes \ $\Ts^{n}$ \ to \ $\Ts^{n+1}$.

A \textit{free simplicial resolution} of a \Pa\ \ $\Gs$ \ 
is defined to be an augmented simplicial \Pa\ \ $\Ad\to\Gs $, \ such
that \

\begin{enumerate}
\renewcommand{\labelenumi}{(\roman{enumi})}
\item \ $\Ad$ \ is a free simplicial \Pa, 
\item in each degree \ $k\geq 1$, \ the homotopy groups of 
the simplicial group \ $(\Ad)_{k}$ \ vanish in dimensions \ $n\geq 1$, \ and 
the augmentation induces an isomorphism \ $\pi _{0} (\Ad)_{k}\cong G_{k}$.
\end{enumerate}
\end{defn}

Such resolutions always exist, for any \Pa\ \ $\Gs$ \ -- \ see
\cite[II, \S 4]{QuH}, or the explicit construction in \cite[\S 4.3]{BlaH}.

\begin{defn}\label{drbs}\stepcounter{subsection}
For any \ $X\in\G$, \ a simplical object \ $\Wd\in s\G$ \ equipped with
an augmentation \ $\varepsilon:W_{0}\to X$ \ is called a 
\textit{resolution of $X$ by spheres\/} 
if each $\W_{n}$ \ is homotopy equivalent to a wedge of $\G$-spheres, 
and \ $\pis\Wd\to\pis X$ \ is a  free simplicial resolution of \Pa s.
\end{defn}

\begin{example}\label{esr}\stepcounter{subsection}
One example of such a resolution by spheres is provided by Stover's 
construction;  we shall need a variant in $\G$ \ (as in \cite[\S 5]{BlaL}), 
rather than the original version of \cite[\S 2]{StoV}, in \ $\Ta$. \ (The 
argument from this point on would actually work equally well in \ $\Ta$; \ 
but we have already chosen to work in $\G$, in order to facilitate the proof of 
Proposition \ref{pone}).

Define a comonad \ $V:\G\to\G$ \ for \ $G\in\G$ \ by 
\setcounter{equation}{\value{thm}}\stepcounter{subsection} 
\begin{equation}\label{ethree}
VG = \ \coprod _{k=0}^{\infty } \ 
               \coprod _{\phi\in \Hom_{\G}(\gS{k},G)} \ \gS{k}_{\phi} \ \ 
               \bigcup \ \ \coprod _{k=0}^{\infty } \ 
               \coprod _{\Phi\in \Hom_{\G} (\gD{k+1},G)} \ \gD{k+1}_{\Phi},
\end{equation}
\setcounter{thm}{\value{equation}}

\noindent where \ $\gD{k+1}_{\Phi}$, \ the $\G$-disc indexed by \ 
$\Phi:\gD{k+1}\to G$, \ is attached to \ $\gS{k}_{\phi}$, \ 
the $\G$-sphere indexed by \ $\phi=\Phi\rest{\partial \gD{k+1}}$,  \ by
identifying \ $\partial \gD{k+1}\DEF F\partial\bD{k}$ \ with \ $\gS{k}$ \ (see \S
\ref{dgs} above).  The coproduct here is just the (dimensionwise) free
product of groups; the counit \ $\varepsilon:VG\to G$ \ of the comonad $V$ is 
``evaluation of indices'', and the comultiplication \ 
$\vartheta:VG\hra V^{2}G$ \ is as in \S \ref{dfpc}.

Now given \ $X\in \G$, \  define \ $\Qd\in s\G$ \ by setting \ 
$Q_{n}=V^{n+1}X$, \ with face and degeneracy maps induced by the counit and 
comultiplication respectively (cf.\ \cite[App., \S 3]{GodT}). \ The counit 
also induces an augmentation \ $\varepsilon:\Qd\to X$; \ and this is in fact 
a resolution of $X$ by spheres (see \cite[Prop.\ 2.6]{StoV}).
\end{example}

\begin{remark}\label{rsr}\stepcounter{subsection}
Note that we need not use the $\G$-sphere and disk \ $\gS{k}$ \ and \ $\gD{k}$ \ 
of \S \ref{dgs} in this construction; we can replace it by any other homotopy 
equivalent cofibrant pair of simplicial groups, so in particular by \ 
$(F\hD{k},F\hS{k-1})$ \ for any pair of simplicial sets \ 
$(\hD{k},\hS{k-1})\simeq(\bD{k},\bS{k-1})$. 
\end{remark}

\subsection{The Quillen spectral sequence}
\label{sqss}\stepcounter{thm}
A resolution by spheres \ $\Wd\to X$ \ is in fact a resolution (i.e., 
cofibrant replacement) for the constant simplicial object \ $\co{X}\in s\G$ \ 
(i.e., \ $c(X)_{n}=X$, \ $d_{i}=s_{j}=id_{X}$) \ in an appropriate model 
category structure on \ $s\G$ \ -- \ see \cite{DKStE} and \cite{BlaAI}. \ 
However, we shall not need this fact; for our purposes it suffices to recall 
that for any bisimplicial group \ $\Wd\in s\G$, \ there is a first quadrant 
spectral sequence with
\setcounter{equation}{\value{thm}}\stepcounter{subsection} 
\begin{equation}\label{efour}
E^{2}_{s,t}=\pi_{s}(\pi_{t}\Wd) \Rightarrow \pi_{s+t}\diag\Wd
\end{equation}
\setcounter{thm}{\value{equation}}

\noindent converging to the diagonal \ $\diag\Wd\in \G$, \ defined \ 
$(\diag\Wd)_{k}=(\W_{k})^{int}_{k}$ \ (see \cite{QuS}). \ Thus if \ 
$\Wd\to X$ \ is a resolution by spheres, the spectral sequence collapses, and 
the natural map \ $\W_{0}\to\diag\Wd$ \ induces an isomorphism \ 
$\pis\X\cong\pis(\diag\Wd)$. \ Combined with the fact that \ $\pis\Wd$ \ is
a resolution (in \ $s\PAlg$) \ of \ $\pis\X$, \ this simple result has many
applications \ -- \ see for example \cite{BlaH}, \cite{DKSStP}, and 
\cite{StoV}.

\begin{defn}\label{dcw}\stepcounter{subsection}
A \textit{CW complex} over a pointed category $\C$ is a simplicial object \ 
$\Rd\in s\C$, \ together with a sequence of objects \ $\bar{R}_{n}$ \ 
($n=0,1,\dotsc$) \ such that \ $R_{n}\cong\bar{R}_{n}\amalg L_{n}\Rd$ \ 
(\S \ref{dmat}), \ and \ $d^{n}_{i}\rest{\bar{R}_{n}}=0$ \ for \ 
$1\leq i\leq n$. \ The objects \ $(\bar{R}_{n})_{n=0}^{\infty}$ \ are called 
a \textit{CW basis} for \ $\Rd$, \ and \ 
$\bar{d}^{n}_{0}\DEF d_{0}\rest{\bar{R}_{n}}$ \ is called the $n$-th attaching
map for \ $\Rd$. \  

One may then describe \ $\Rd$ \ explicitly in terms of its CW basis by  
\setcounter{equation}{\value{thm}}\stepcounter{subsection} 
\begin{equation}\label{efive}
R_{n}\cong\coprod_{0\leq \lambda \leq n}\ \coprod _{I \in\I_{\lambda ,n}}\ 
\bar{R}_{n-\lambda }
\end{equation}
\setcounter{thm}{\value{equation}}

\noindent where \ $\I_{\lambda ,n}$ \ is the set of sequences $I$ 
of $\lambda$ non-negative integers \ $i_{1} < i_{2} < \dotsc < i_{\lambda}$ \ 
($<n$), \ with \ 
$s_{I}=s_{i_{\lambda}}\circ\dotsb\circ s_{i_{0}}$ \ the corresponding
$\lambda$-fold degeneracy (if \ $\lambda=0$, \ $s_{I}=id$). \ See
\cite[5.2.1]{BlaD} and \cite[p.\ 95(i)]{MayS}.
\end{defn}

Such CW bases are convenient to work with in many situations; but they are most 
useful when each basis object is \emph{free}, in an appropriate sense. In 
particular, if \ $\C=\PAlg$, \ we have the following

\begin{defn}\label{dcwr}\stepcounter{subsection}
A \textit {CW resolution} of a \Pa\ \ $\Gs$ \ is a CW complex \ 
$\Ad\in s\PAlg$, \ with CW basis \ $(\bar{A}_{n})_{n=0}^{\infty}$ \ and 
attaching maps \ $\bar{d}^{n}_{0}:\bar{A}_{n}\to Z_{n-1}\Ad$, \ such that each \ 
$\bar{A}_{n}$ \ is a free \Pa, and each attaching map
$d^{n}_{0}\rest{C_{n}\Ad}$ \ is onto \ $Z_{n-1}\Ad$ \ (for \ $n\geq 0$, \ 
where we let \ $\bar{d}^{0}_{0}$ \ denote the augmentation \ 
$\varepsilon:\Ad\to \Gs$ \ and \ $Z_{-1}\Ad:=\Gs$). \ Compare \cite[\S 5]{BlaD}.
\end{defn}

Every \Pa\ has a CW resolution (\S \ref{dcwr}), as was shown in 
\cite[4.4]{BlaH}: \ for example, one could take the graded set of generators \ 
$\bar{T}^{n}_{\ast}$ \ for \ $\bar{A}_{n}$ \ to be equal to the graded 
\textit{set} \ $\pis Z_{n-1}\Ad$.

\begin{defn}\label{dcwg}\stepcounter{subsection}
$\Qd\in s\G$ \ is called a \textit{CW resolution by spheres} of \ $X\in\G$ \ 
if \ $\Qd\to X$ \ is a resolution by spheres (Def.\ \ref{drbs}), and \ 
$\Qd$ \ is a CW complex with CW basis \ $(\bar{Q}_{n})_{n=0}^{\infty}$), \ 
such that each \ $\bar{Q}_{n}\in\F$ \ (i.e., \ $\bar{Q}_{n}$ \ is homotopy 
equivalent to a wedge of spheres). The concept is defined analogously for 
$X\in\Ss$ \ or \ $X\in\Ta$.
\end{defn}

\begin{remark}\label{rpa}\stepcounter{subsection}\
Closely related to the problem of realizing abstract \Pa s (\S \ref{dpa})
is that of realizing a free simplicial \Pa\ \ $\Ad\in s\PAlg$: \ this is 
because, as noted in \S \ref{dfsr}, every \ $\Gs\in\PAlg$ \ has a free 
simplicial resolution \ $\Ad\to\Gs$; \ if it can be realized by a simplicial
space \ $\Wd\in s\Tc$ \ -- \ or equivalently, via \eqref{ezero}, by a 
bisimplicial space or group \ -- \ then the spectral sequence \eqref{efour} 
implies that \ $\pis\diag\Wd\cong\Gs$. \ However, not every \Pa\ is realizable
(see \cite[\S 8]{BlaHH} or \cite[Prop.\ 4.3.6]{BlaO}).

It would nevertheless be very useful to know the converse: namely, that any 
free resolution of a \textit{realizable} \Pa\ is itself realizable. This was 
mistakenly quoted as a theorem in \cite[\S 6]{BlaHH}, where it was needed to 
make the obstruction theory for realizing \Pa s described there of any 
practical use \ -- \ and appeared as a conjecture in \cite[\S 4]{BlaHO}, in 
the context of an obstruction theory for a space to be an $H$-space. 
\end{remark}

In order to show that this conjecture is in fact true, we need several
preliminary results:

%
%
\begin{prop}\label{ptwo}\stepcounter{subsection}
Every CW resolution \ $\Ad\to\pis X$ \ of a \textit{realizable} \Pa\ 
embeds in \ $\pis\Qd$ \ for some resolution by spheres \ $\Qd\to X$.
\end{prop}

\begin{proof}
To simplify the notation, we work here with topological spaces, rather
than simplicial groups, changing back to $\G$ if necessary via the adjoint 
pairs of \S \ref{ssg}.

Given a free simplicial \Pa\ resolution \ $\Ad\to\Js$ \ with CW basis \ 
$(\bar{A}_{n})_{n=0}^{\infty}$, \ where \ $\Js=\pis\X$ \ for some \ 
$\X\in\Ta$, \ and \ $\bar{A}_{n}$ \ is the free \Pa\ generated by the graded
set \ $\Ts^{n}$, \ let $\mu$ denote the cardinality of \ 
$\coprod_{n=0}^{\infty}\coprod_{k=0}^{\infty} T^{n}_{k}$, \ and set \ 
$\X'\DEF \X\vee\bigvee_{n=0}^{\infty}\bigvee_{\lambda<\mu}\bD{n}$. \ 
Define new ``spheres'' and ``disks'' of the form \ 
$\hS{n}\DEF\bS{n}\vee\bigvee_{n=0}^{\infty}\bigvee_{\lambda<\mu}\bD{n}$ \ 
and \ $\hD{n}\DEF\hS{n}\vee\bD{n}$. \ 
(This is to ensure that there will be at least $\mu$ different representatives 
for each homotopy class in \ $\pis\X'$ \ or \ $\pis\hS{n}$.)

By remark \S \ref{rsr} above, if we use the construction of \S \ref{esr} in \ 
$\Ta$ \ (or in $\G$, \textit{mutatis mutandis}) with these ``spheres'' and 
``disks'', and apply it to the space \ $\X'$, \ rather than to $\X$, we obtain a 
resolution by spheres \ $\Qd\to\X'$.

We define \ $\phi:\Ad\hra \pis\Qd$ \ by induction on the simplicial dimension;
it suffices to produce for each \ $n\geq 0$ \ an embedding \ 
$\bar{\phi}_{n}:\bar{A}_{n}\hra C_{n}\pis\Qd$ \ commuting with \ $d_{0}$. \ 
If we denote \ $\varepsilon^{A}:A_{0}\to\pis\X\cong\pis\X'$ \ by \ 
$\bar{d}_{0}^{0}:C_{0}\Ad\to Z_{-1}\Ad=:A_{-1}$ \ and set \ 
$\phi_{-1}=id_{\pis X}$, \ then we may assume by induction we have a 
monomorphism \ $\phi_{n-1}:A_{n-1}\hra \pis Q_{n-1}$ \ (taking generators to
generators, and commuting with face and degeneracy maps).

For each \Pa\ generator \ $\iota_{\alpha}$ \ in \ $(\bar{A}_{n})_{k}$, \ if \ 
$d_{0}(\iota_{\alpha})\neq 0$ \ then \ 
$\phi_{n-1}(d_{0}(\iota_{\alpha}))\in Z_{n-1}\pi_{k}\Qd$ \ is represented by 
some \ $g:\hS{k}\to Q_{n-1}$, \ and we can choose distinct (though perhaps
homotopic) maps $g$ for different generators \ $\iota_{\alpha}$ \ by our
choice of \ $\hS{k}$. \ Then by \eqref{ethree} there is a wedge summand \ 
$\hS{k}_{g}$ \ in \ $Q_{n}=VQ_{n-1}$ \ (with no disks attached), \ and the
corresponding free \Pa\ coproduct summand \ $\pis\hS{k}_{g}$ \ in \ 
$\pis Q_{n}$, \ generated by \ $\iota_{g}$, \ has \ 
$d_{0}(\iota_{g})=[g]\in\pi_{k} Q_{n-1}$ \ and \ 
$d_{i}(\iota_{g})=\iota_{d_{i-1}g}=0\in\pi_{k} Q_{n-1}$ \ for \ 
$1\leq i\leq n$ \ by \S \ref{esr}, since \ 
$[g]=\phi_{n-1}(d_{0}(\iota_{\alpha}))\in Z_{n-1}\pi_{k}\Qd$ \ and thus \ 
$d_{i}[g]=[d_{i}g]=0$, \ and spheres indexed by nullhomotopic maps have 
disks attached to them. \ We see that \ $\iota_{g}\in C_{n}\pi_{k}\Qd$, \ 
so we may define \ $\bar{\phi}_{n}(\iota_{\alpha})=\iota_{g}$.

If \ $d_{0}(\iota_{\alpha})=0$, \ then all we need are enough distinct
\Pa\ generators in \ $Z_{n}\pis\Qd$: \ we cannot simply take \ 
$\iota_{g}$ \ for nullhomotopic \ $g:\bS{k}\to Q_{n-1}$, \ because of the
attached disks; but we can proceed as follows:

Since \ $\hD{k}=C\hS{k}\vee\bD{k}$ \ and \ 
$\X'=\X\vee\bigvee_{i=0}^{\infty}\bigvee_{\lambda<\mu}\bD{i}$, \ we have
$\mu$ distinct nonzero maps \ $F_{\lambda}:\hD{k}\to\X'$ \ with \ 
$F_{\lambda}\rest{C\hS{k}}=\ast$. \ Define \ $H_{+}=F_{\lambda}$, \ 
$H_{-}=\ast$; \ then
$\bS{k}_{H}\DEF\hD{k}_{H^{+}}\cup_{\hS{k-1}_{\ast}}\hD{k}_{H^{-}}$ \ is, up to 
homotopy, a sphere wedge summand in \ $Q_{0}$, \ and thus \ 
$\iota_{H_{\lambda}}\in\pi_{k} Q_{0}$ \ is a \Pa\ generator mapping to $0$ 
under the augmentation. Similalry, define
$\bS{k}_{G_{\lambda}}\DEF\hD{k}_{G^{+}}\cup_{\hS{k-1}_{\ast}}\hD{k}_{G^{-}}$ \ 
in \ $Q_{1}$ \ by \ $G^{+}=\ast$, \ $G_{-}=\ast\bot\iota^{k}$ \ 
where \ $\iota^{k}$ \ is a homoeomorphism onto the summand \ $\bD{k}$ \ in \ 
$\hD{k}_{H^{-}_{\lambda}}$. \ Then \ $G_{\lambda}\sim\ast$ \ and \ 
$G_{\lambda}\neq \ast$ \ but \ $H\circ G=\ast$; \ thus \ 
$\iota_{H_{\lambda}}$ \ is a \Pa\ generator in \ $Z_{1}\pi_{k}\Qd$. \ 
By thus alternating the $+$ and $-$ we produce $\mu$ distinct \Pa\ generators 
in \ $Z_{n}\pis\Qd$ \ for each $n$.
\end{proof}

\begin{remark}\label{rref}\stepcounter{subsection}
The referee has suggested an alternative proof of this Proposition, which may be 
easier to follow: rather than ``fattening'' the spheres and disks, we can modify 
the Stover construction of \eqref{ethree} by using $\mu$ copies of each sphere or 
disk for each \ $\phi\in \Hom_{\G}(\gS{k},G)$ \ or \ 
$\Phi\in \Hom_{\G} (\gD{k+1},G)$, \ respectively. The proof of 
\cite[Prop.\ 2.6]{StoV} still goes through, and so does the argument for 
embedding \ $\Ad$ \ in \ $\pis\Qd$ \ above.
\end{remark}

%
%
\begin{prop}\label{pthree}\stepcounter{subsection}
Any free simplicial \Pa\ \ $\Ad$ \ has a (free) CW basis \ 
$(\bar{A}_{n})_{n=0}^{\infty}$.
\end{prop}

\begin{proof}
Start with \ $\bar{A}_{0}=A_{0}$.
For \ $n\geq 1$, \ assume \  
$A_{n}=\coprod_{k=0}^{\infty}\coprod_{\tau\in T^{n}_{k}}\pis\gS{k}$. \ 
By Definition \ref{dfsr}, \ 
$\Ts^{n}\cong \bar{T}^{n}_{\ast}\ \cup \bigcup_{0\leq \lambda \leq n}\ 
\bigcup_{I \in\I_{\lambda ,n}}\ \hat{T}^{n-\lambda}_{\ast}$ \ 
(as in \S \ref{efive}), \ so we can set \ $\hat{A}_{n}=
\coprod_{k=0}^{\infty}\coprod_{\tau\in \hat{T}^{n}_{k}}\pis\gS{k}$; \ 
but \ $d_{i}\rest{\hat{A}_{n}}$ \ need not vanish for \ $i\geq 1$.

However, given \ $\tau\in \hat{T}^{n}_{k}$, \ we may define \ 
$\tau_{i}\in (A_{n})_{k}^{int}$ \ inductively, starting with \ 
$\tau_{0}=\tau$, \ by \ $\tau_{i+1}=\tau_{i}s_{n-i-1}d_{n-i}\tau_{i}^{-1}$ \ 
(face and degeneracy maps taken in the external direction); \ 
we find that \ $\bar{\tau}\DEF\tau_{n}$ \ is in \ $C_{n}\Ad$. \ If we define \ 
$\bar{\varphi}:\hat{T}^{n}_{\ast}\to A_{n}$ \ by \ 
$\varphi(\tau)=\bar{\tau}$, \ by the universal property of free \Pa s this 
extends to a map \ $\varphi:\hat{A}_{n}\to A_{n}$, \ which together with
the inclusion \ $\sigma_{n}:L_{n}\Ad\hra A_{n}$ \ yields a map \ 
$\psi:A_{n}\to A_{n}$ \ which is an isomorphism by the Hurewicz Theorem
(cf.\ \cite[Lemma 2.5]{BlaL}). \ Thus we may set \ 
$\bar{A}_{n}\DEF \varphi(\hat{A}_{n})$, \ that is, the free \Pa\ generated 
by \ $\{\bar{\tau}\}_{\tau\in\hat{T}^{n}_{\ast}}$. \ 
Compare \cite[\S 3]{KanR}.
\end{proof}

%
%
\begin{thm}\label{tone}\stepcounter{subsection}
Every free simplicial \Pa\ resolution \ $\Ad\to\pis X$ \ of a realizable 
\Pa\ \ $\pis X$ \ is itself realizable by a CW resolution \ $\Rd\to X$ \ 
in \ $s\G$.
\end{thm}

\begin{proof}
By Propositions \ref{ptwo} and \ref{pthree} we may assume \ $\Ad$ \ has a (free)
CW basis \ $(\bar{A}_{n})_{n=0}^{\infty}$, \ and that there is a resolution
by spheres \ $\Qd\to X$ \ (in \ $s\G$) \ and an embedding of simplicial 
\Pa s \ $\phi:\Ad\to\Qd$. \ We may also assume that \ $\Qd$ is fibrant 
(\S \ref{dfib}), with \ $\varepsilon^{Q}:Q_{0}\to X$ \ a fibration. 
We shall actually realize $\phi$ \ by a map of bisimplicial groups \ 
$f:\Rd\to\Qd$.

Note that once \ $\Rd$ \ has been defined through
simplicial dimension $n$, for any \ $k\geq 0$ \ we have a commutative diagram

%
%
\begin{picture}(400,90)(-5,0)
%
%
\put(0,60){$\pi_{k}C_{n}\Rd$}
\put(43,65){\vector(1,0){34}}
\put(48,70){$(d_{0})_{\#}$}
\put(80,60){$\pi_{k}Z_{n-1}\Rd$}
\put(15,54){\vector(0,-1){35}}
\put(2,37){$\iota_{\star}$}
\put(18,37){$\cong$}
\put(95,54){\vector(0,-1){35}}
\put(72,37){$\rho_{n-1}$}
\put(0,10){$C_{n}\pi_{k}\Rd$}
\put(43,15){\vector(1,0){34}}
\put(55,20){$d_{0}^{n}$}
\put(80,10){$Z_{n-1}\pi_{k}\Rd$}
%
%
\put(132,65){\vector(1,0){43}}
\put(136,70){$(j_{n-1})_{\#}$}
\put(180,60){$\pi_{k}C_{n-1}\Rd$}
\put(195,54){\vector(0,-1){35}}
\put(182,37){$\iota_{\star}$}
\put(199,37){$\cong$}
\put(132,15){\vector(1,0){43}}
\put(140,20){inc.}
\put(180,10){$C_{n-1}\pi_{k}\Rd$}
%
%
\put(233,65){\vector(1,0){33}}
\put(238,70){$(d_{0})_{\#}$}
\put(270,60){$\pi_{k}Z_{n-2}\Rd$}
\put(290,54){\vector(0,-1){35}}
\put(267,37){$\rho_{n-2}$}
\put(233,15){\vector(1,0){33}}
\put(245,20){$d_{0}^{n-1}$}
\put(270,10){$Z_{n-2}\pi_{k}\Rd$}
%
%
\put(323,65){\vector(1,0){43}}
\put(328,70){$(j_{n-2})_{\#}$}
\put(370,60){$\pi_{k}C_{n-2}\Rd$}
\put(385,54){\vector(0,-1){35}}
\put(372,37){$\iota_{\star}$}
\put(389,37){$\cong$}
\put(323,15){\vector(1,0){43}}
\put(330,20){inc.}
\put(370,10){$C_{n-2}\pi_{k}\Rd$}
\end{picture}

\noindent (obtained by fitting together three of the long exact sequences
of the fibrations \eqref{etwo}). The vertical maps are induced by the 
inclusions \ $C_{n}\Rd\hra R_{n}$, \ and so on  \ -- \ see Proposition 
\ref{pone}.  

The only difficulty in constructing \ $\Rd$ \ is that Proposition
\ref{pone} does not hold for \ $Z_{n}$ \ -- \ i.e., the maps \ $\rho_{n}$ \ in
the above diagram in general need not be isomorphisms \ -- \ so we may have an
element in \ $Z_{n}\Ad$ \ represented by \ 
$\alpha\in C_{n}\pis\Rd=\pis C_{n}\Rd$ \ with \ 
$(d^{n}_{0})_{\#}(\alpha)\neq 0$ \ (but of course \ 
$(j_{n-1})_{\#}(d^{n}_{0})_{\#}(\alpha)=0$). \ In this case we could not have \ 
$\beta\in \pis C_{n+1}\Rd=C_{n+1}\Ad$ \ with \ 
$(j_{n})_{\#}(d^{n+1}_{0})_{\#}(\beta)=\alpha$, \ so \ $\pis\Rd$ \ would not be 
acyclic.

It is in order to avoid this difficulty that we need the embedding $\phi$,
since by definition this cannot happen for \ $\Qd$: \ we know that \ 
$d_{0}^{n}:C_{n}\pis\Qd\to Z_{n-1}\pis\Qd$ \ is surjective for each \ $n>0$, \
so \ $\rho_{n-1}:\pis Z_{n-1}\Qd\to Z_{n-1}\pis\Qd$ \ is, too, which implies
that for each \ $n>0$:
\setcounter{equation}{\value{thm}}\stepcounter{subsection} 
\begin{equation}\label{esix}
\Image\{(d_{0}^{n+1})_{\#}:\pis C_{n+1}\Qd\to \pis Z_{n}\Qd\}\cap
\Ker\{(j_{n})_{\#}:\pis Z_{n}\Qd\to \pis C_{n}\Qd\}=0
\end{equation}
\setcounter{thm}{\value{equation}}

\noindent which we shall call \textit{Property \eqref{esix} for \ 
$Z_{n}\Qd$}. (This implies in particular that \ 
$Z_{n}\pis\Qd=\Ker\{(d_{0}^{n})_{\#}:\pis C_{n}\Qd\to Z_{n-1}\Qd\}$.)

Note that given any fibrant \ $\Kd\in s\G$ \ having Property \eqref{esix} 
for \ $Z_{m}\Kd$ \ for each \ $0<m\leq n$, \ if we consider the long exact 
sequence of the fibration \ $d^{m}_{0}:C_{m}\Kd\to Z_{m-1}\Kd$:
%
\setcounter{equation}{\value{thm}}\stepcounter{subsection} 
\begin{equation}\label{eseven}
\dotsc \pi_{k+1}C_{m}\Kd \xra{(d^{m}_{0})_{\#}} \pi_{k+1}Z_{m-1}\Kd
\xra{\partial^{m-1}} \pi_{k}Z_{m}\Kd \xra{(j_{m-1})_{\#}} \pi_{k}C_{m-1}\Kd 
\dotsc,
\end{equation}
\setcounter{thm}{\value{equation}}
we may deduce that \ 
%
\setcounter{equation}{\value{thm}}\stepcounter{subsection} 
\begin{equation}\label{eeight}
\partial^{m}\rest{\Image(\partial^{m-1}) } \ \ \ 
\text{is one-to-one, and surjects onto} \ \Image(\partial^{m}) 
\end{equation}
\setcounter{thm}{\value{equation}}
for \ $0<m\leq n$\vsm .

We now construct \ $\Rd$ \ by induction on the simplicial dimension\vsm :
\begin{enumerate}
\renewcommand{\labelenumi}{(\roman{enumi})}
\item First, choose a fibration \ $\varepsilon^{R}:R_{0}\to X$ \ realizing \
$\varepsilon^{A}:A_{0}\to \pis X$. \ By \S \ref{dfpa}, there is a map \ 
$f_{0}':R'_{0}\to Q_{0}$ \ realizing \ $\phi_{0}$,  \ so \ 
$\varepsilon^{Q}\circ f_{0}'\sim \varepsilon^{R}$; \ since \ 
$\varepsilon^{Q}$ \ is a fibration, we can change \ $f'_{0}$ \ to \ 
$f_{0}:R_{0}\to Q_{0}$ \ with \ 
$\varepsilon^{Q}\circ f_{0}'=\varepsilon^{R}$\vsm .

\item Let \ $Z_{0}\Rd$ \ denote the fiber of \ $\varepsilon^{R}$.
Since \ $\varepsilon^{R}_{\#}=\varepsilon^{A}$ \ is a surjection, we have \ 
$\pis Z_{0}\Rd=\Ker(\varepsilon^{R}_{\#})=Z_{0}\Ad$, \ and \ $d^{A}_{0}$ \ maps \ 
$C_{1}\Ad$ \ onto \ $Z_{0}\Ad$, \ so  \ 
$\bar{d}_{0}^{A}:\bar{A}_{1}\to A_{0}$ \ factors through \ $\pis Z_{0}\Rd$, \ 
and we can thus realize it by a map \ 
$\bar{d}_{0}^{R}:\bar{R}_{1}\to Z_{0}\Rd$. \ Set \ 
$R_{1}'\DEF \bar{R}_{1}\amalg L_{1}\Rd$ \ (so \ $\pis R_{1}'\cong A_{1}$), \ 
with \ $\delta'_{1}:R_{1}'\to M_{1}\Rd=R_{0}\times R_{0}$ \ equal to \ 
$(\bar{d}^{R}_{0},0)\bot \Delta$, \ and change \ $\delta'_{1}$ \ to a 
fibration \ $\delta_{1}:R_{1}\to M_{1}\Rd$. \ Again we can realize \ 
$\phi_{1}:A_{1}\to\pis Q_{1}$ \ by \ $f_{1}:R_{1}\to Q_{1}$ \ with \ 
$\delta_{1}^{Q}\circ f_{1}=f_{0}\circ\delta_{1}^{R}$, \ since \ 
$\delta_{1}^{Q}$ \ is a fibration; so we have defined \ 
$\tau_{1}f:\tau_{1}\Rd\to\tau_{1}\Qd$ \ realizing \ $\tau_{1}\phi$\vsm .

\item Now assume we have \ $\tau_{n}f:\tau_{n}\Rd\to\tau_{n}\Qd$ \ 
realizing \ $\tau_{n}\phi$, \ with Property \eqref{esix} holding for \ 
$Z_{m}\Rd$ \ for \ $0<m<n$.

For each \Pa\ generator \ $\alpha\in\bar{A}_{n+1}$ \ (in degree $k$, say),
\eqref{esix} implies that \ 
$d_{0}^{n+1}(\alpha)\in \Ker(d_{0}^{n})=\Ker((d_{0}^{R_{n}})_{\#})\subset 
(C_{n}\Ad)_{k}=\pi_{k}C_{n}\Rd$, \ so by the exactness of \eqref{eseven} we 
can choose \ $\beta\in\pi_{k}\Z_{n}\Rd$ \ such that \ 
$(j_{n})_{\#}\beta=d_{0}^{n+1}(\alpha)$. \ This allows us to define \ 
$\bar{d}_{0}^{R}:\bar{R}_{n+1}\to Z_{n}\Rd$ \ so that \ 
$(j_{n})_{\#}(\bar{d}_{0}^{R})_{\#}$ \ realizes \
$(\text{inc.})\circ\bar{d}_{0}^{A}:\bar{A}_{n+1}\to C_{n}\Ad$, \ as well as \ 
$\bar{f}_{n+1}:\bar{R}_{n+1}\to C_{n}\Qd$ \ realizing \ 
$\phi_{n+1}\rest{\bar{A}_{n+1}}$. \ Because \ 
$\bar{A}_{n+1}=\pis\bar{R}_{n+1}$ \ is a free \Pa, this implies the 
homotopy-commutativity of the outer rectangle in

%
%
\begin{picture}(120,135)(-120,-10)
%
%
\put(0,110){$\bar{R}_{n+1}$}
\put(27,115){\vector(1,0){48}}
\put(45,120){$\bar{f}_{n+1}$}
\put(80,110){$C_{n+1}\Qd$}
\put(15,103){\vector(0,-1){42}}
\put(0,85){$\bar{d}_{0}^{R}$}
\put(95,103){\vector(0,-1){42}}
\put(98,85){$d_{0}^{Q}$}
\put(0,50){$Z_{n}\Rd$}
\multiput(30,55)(3,0){15}{\circle*{.5}}
\put(75,55){\vector(1,0){3}}
\put(45,60){$Z_{n}f$}
\put(80,50){$Z_{n}\Qd$}
\put(15,45){\vector(0,-1){35}}
\put(0,27){$j_{n}^{R}$}
\put(95,45){\vector(0,-1){35}}
\put(98,27){$j_{n}^{Q}$}
\put(0,0){$C_{n}\Rd$}
\put(30,5){\vector(1,0){45}}
\put(50,10){$C_{n}f$}
\put(80,0){$C_{n}\Qd$}
\end{picture}

\noindent (as well as the lower square, by the induction hypothesis).
Thus \ 
$j_{n}^{Q}\circ Z_{n}f\circ \bar{d}_{0}^{R}\sim 
j_{n}^{Q}\circ d_{0}^{Q}\circ \bar{f}_{n+1}$, \ so \ 
$(j_{n}^{Q})_{\#}\circ (Z_{n}f)_{\#}\circ (\bar{d}_{0}^{R})_{\#}=
(j_{n}^{Q})_{\#}\circ(d_{0}^{Q})_{\#}\circ(\bar{f}_{n+1})_{\#}$. \ By
\eqref{esix} this implies \ 
$(Z_{n}f)_{\#}\circ (\bar{d}_{0}^{R})_{\#}=
(d_{0}^{Q})_{\#}\circ(\bar{f}_{n+1})_{\#}$, \ so (since \ 
$\pis\bar{R}_{n+1}$ \ is a free \Pa) \ also \ 
$Z_{n}f\circ \bar{d}_{0}^{R}\sim d_{0}^{Q}\circ \bar{f}_{n+1}$ \ -- \ which
means that we can choose \ $\bar{f}_{n+1}$ \ so that \ 
$Z_{n}f\circ \bar{d}_{0}^{R}=d_{0}^{Q}\circ \bar{f}_{n+1}$ \ (since \ 
$d_{0}^{Q}$ \ is a fibration). Thus if we set \ 
$\bar{\delta}^{R}_{n+1}:\bar{R}_{n+1}\to M_{n+1}\Rd$ \ to be \ 
$(\bar{d}_{0}^{R},0,\dotsc,0)$, \ we have \ 
$M_{n+1}f\circ\bar{\delta}_{n+1}^{R}=\delta_{n+1}^{Q}\circ\bar{f}_{n+1}$. \ 

If \ $\psi^{R}_{n+1}\DEF\delta^{R}_{n+1}\circ\sigma^{R}_{n+1}$ \ (in the 
notation of \S \ref{dmat} \& \ref{dlat}) \ we set \ 
$R'_{n+1}\DEF \bar{R}_{n+1}\amalg L_{n+1}\Rd$, \ and define \ 
$\delta'_{n+1}:R'_{n+1}\to M_{n+1}\Rd$, \ and \ $f'_{n+1}:R'_{n+1}\to Q_{n+1}$ \ 
respecyively by \ 
$\delta'_{n+1}\DEF (\bar{\delta}^{R}_{n+1}\bot\psi^{R}_{n+1})$ \ and \ 
$f'_{n+1}\DEF (\bar{f}_{n+1}\bot L_{n+1}f)$. \ We see that \ 
$(f'_{n+1})_{\#}=\phi_{n+1}$ \ and \ 
$M_{n+1}f\circ\delta'_{n+1}=\delta_{n+1}^{Q}\circ f'_{n+1}$, \ and this will
still hold if we change \ $\delta'_{n+1}$ \ into a fibration, and extend \ 
$f'_{n+1}$ \ to \ $f_{n+1}:R_{n+1}\to Q_{n+1}$. \ This defines \ 
$\tau_{n+1}f:\tau_{n+1}\Rd\to\tau_{n+1}\Qd$ \ realizing \ $\tau_{n+1}\phi$\vsm .
\item It remains to verify that \ $\tau_{n+1}\Rd$ \ so defined satisfies 
\eqref{esix}. However, \eqref{eeight} implies that we have a map of short exact 
sequences:

%
%
\begin{picture}(350,95)(-10,-10)
%
%
\put(0,60){$0$}
\put(10,65){\vector(1,0){22}}
\put(35,60){$\Image(\partial^{n-1}_{R})$}
\put(87,65){\vector(1,0){43}}
\put(100,70){inc.}
\put(135,60){$\pi_{k}Z_{n}\Rd$}
\put(55,53){\vector(0,-1){39}}
\put(57,32){$f_{\ast}$}
\put(155,53){\vector(0,-1){39}}
\put(158,32){$(Z_{n}f)_{\#}$}
\put(0,0){$0$}
\put(10,5){\vector(1,0){22}}
\put(35,0){$\Image(\partial^{n-1}_{Q})$}
\put(87,5){\vector(1,0){45}}
\put(100,10){inc.}
\put(135,0){$\pi_{k}Z_{n}\Qd$}
%
%
\put(177,65){\vector(1,0){45}}
\put(225,60){$\Image((j_{n}^{R})_{\#})\cong Z_{n}\Ad$}
\put(320,65){\vector(1,0){30}}
\put(355,60){$0$}
\put(177,5){\vector(1,0){45}}
\put(225,0){$\Image((j_{n}^{Q})_{\#})\cong Z_{n}\pi_{k}\Qd$}
\put(335,5){\vector(1,0){15}}
\put(355,0){$0$}
\put(295,53){\vector(0,-1){39}}
\put(298,32){$Z_{n}\phi$}
\end{picture}

\noindent in which the left vertical map is an isomorphism and the right map 
is one-to-one, so \ $(Z_{n}f)_{\#}$ \ is one-to-one, too. Therefore, 
$\Ker((j_{n}^{R})_{\#})=\Ker((j_{n}^{R})_{\#})\cap \pis Z_{n}\Rd$, \ which
implies that Property \eqref{esix} holds for \ $Z_{n}\Rd$, \ too.
\end{enumerate}
This completes the inductive construction of \ $\Rd$.
\end{proof}

We also have an analogous result for maps:
%
%
\begin{thm}\label{ttwo}\stepcounter{subsection}
If \ $\Kd\xra{\varepsilon^{K}}\pis X$ \ and \ 
$\Ld\xra{\varepsilon^{L}}\pis Y$ \ are two free simplicial \Pa\ 
resolutions, \ $g:X\to Y$ \ is a map in $\G$, \ and \ 
$\varphi:\Kd\to\Ld$ \ is a morphism of simplicial \Pa s such that \ 
$\varepsilon^{L}\circ \varphi_{0}=\pis g\circ \varepsilon^{K}$, \ then 
$\varphi$ is realizable by a map \ $f:\Ad\to\Bd$ \ in \ $s\G$.
\end{thm}

\begin{proof}
Choose free CW bases for \ $\Kd$ \ and \ $\Ld$, \ and realize the resulting CW
resolutions by \ $\Ad$ \ and \ $\Bd$ \ respectively, where (as in the proof of
Theorem \ref{tone}) we may assume \ $d_{0}:C_{n}\Bd\to Z_{n-1}\Bd$ \ is a 
fibration for each \ $n\geq 0$. \ $f_{n}:A_{n}\to B_{n}$ \ will be defined by 
induction on $n$: \ $\varphi_{0}:K_{0}\to L_{0}$ \ may be realized by a map \ 
$f'_{0}:A_{0}\to B_{0}$ \ (\S \ref{dfpa}), \ and since \ $\varepsilon^{B}$ \ 
is a fibration and \ 
$\varepsilon^{B}\circ f'_{0}\sim g\circ\varepsilon^{A}$, \ we can choose a 
realization \ $f_{0}$ \ for \ $\varphi_{0}$ \ such that \ 
$\varepsilon^{B}\circ f_{0}=g\circ\varepsilon^{A}$.

In general, \ 
$\bar{\varphi}_{n}=\varphi_{n}\rest{\bar{K}_{n}}:\bar{K}_{n}\to C_{n}\Ld$ \ 
may be realized by a map \ $\bar{f}_{n}:\bar{A}_{n}\to C_{n}\Bd$ \ 
(Proposition \ref{pone}), \ and since \ $d_{0}:C_{n}\Bd\to Z_{n-1}\Bd$ \ is a 
fibration, we may choose \ $\bar{f}_{n}$ \  so \ 
$d_{0}\circ \bar{f}_{n}= Z_{n-1}f\circ d_{0}:\bar{A}_{n}\to Z_{n-1}\Bd$. \ 
By induction this yields a map \ 
$f_{n}=L_{n}f\bot\bar{f}_{n}:A_{n}=L_{n}\Ad\amalg \bar{A}_{n}\to 
L_{n}\Bd\amalg \bar{B}_{n}=B_{n}$ \ such that \ 
$\delta_{n}^{B}\circ f_{n}=M_{n}f\circ \delta_{n}^{A}:A_{n}\to M_{n}\Bd$, \ 
so $f$ is indeed a simplicial morphism (realizing $\phi$).
\end{proof}
%
%
\section{The simplicial bar construction}
\label{csbc}

As an application of Theorem \ref{tone}, we describe an obstruction theory
for determining whether a given space $X$ is, up to homotopy, a
loop space (and thus a topological group \ -- \ see \cite[\S 3]{MilnC1}).
In the next two sections we no longer need to work with simplicial groups, 
so we revert to the more familiar category of topological spaces; we can 
still utilize the results of the previous section via the adjoint pairs of 
\eqref{ezero}.

\begin{defn}\label{ddcs}\stepcounter{subsection}
A $\Delta$-\textit{cosimplicial object\/} \ $E^{\bullet}_{\Delta}$ \
over a category $\C$ is a sequence of objects \ $E^{0},E^{1},\dotsc$, \ 
together with \textit{coface maps\/} \ $d^{i}:E^{n}\to E^{n+1}$ \ for \
$1\leq 1\leq n$ \ satisfying \ $d^{j}d^{i}=d^{i}d^{j-1}$ \ for \ $i<j$ \
(cf.\ \cite{RSaD1}). \  
Given an ordinary cosimplicial object \ $E^{\bullet}$ \ 
(cf.\ \cite[X, 2.1]{BKaH}), \ we let \ $E^{\bullet}_{\Delta}$ \ denote the 
underlying $\Delta$-cosimplicial object (obtained by forgetting the 
codegeneracies).
\end{defn}

\subsection{The cosimplicial James construction}
\label{scsc}\stepcounter{thm}

Given a space \ $\X\in\Ta$, \ we define a $\Delta$-cosimplicial 
space \ $\Udd=U(\X)^{\bullet}_{\Delta}$ \ 
by setting \ $\U^{n}=\X^{n+1}$ \ (the Cartesian product), \ and \ 
$d^{i}(x_{0},\dotsc,x_{n})=(x_{0},\dotsc,x_{i-1},\ast,x_{i},\dotsc,x_{n})$. \ 
Note that \ $\colim\U(X)^{\bullet}_{\Delta}\cong J\X$ \ 
(the James reduced product construction), and 

\begin{fact}\label{fahs}\stepcounter{subsection}
If \ $\lra{\X,m}$ \ is a (strictly) associative $H$-space, we can 
extend \ $\Udd$ \ to a full cosimplicial space \ $\Ud$ \ by setting \ 
$s^{j}(x_{0},\dotsc,x_{n})=(x_{0},\dotsc,m(x_{j},x_{j-1}),\dotsc,x_{n})$\vsm.
\end{fact}

\begin{defn}\label{dedd}\stepcounter{subsection}
Let \ $\Ad$ \ be a $CW$-resolution of the \Pa\ \ $\pis\X=\pis\U^{0}$, \ as in
\S \ref{dcwr}. We construct a $\Delta$-cosimplicial
augmented simplicial \Pa\ \ $\Edd\to\pis\Udd$, \ such that each \ 
$E_{\bullet}^{n}$ \ is a $CW$-resolution of \ 
$\pis\U^{n}=\pis(\X^{n+1})$, \ with $CW$-basis \ 
$\{\bar{E}_{r}^{n}\}_{r=0}^{\infty}$. \ We start by setting \ 
$\bar{E}_{r}^{0}=\bar{C}_{r}^{0}=\bar{A}_{r}$ \ for all \
$r\geq 0$, \ and then define \ $\bar{E}_{r}^{n}$ \ by a double
induction (on \ $r\geq 0$ \ and then on \ $n\geq 0$) \ as

\setcounter{equation}{\value{thm}}\stepcounter{subsection}
\begin{equation}\label{enine}
\bar{E}_{r}^{n}~ =~ \coprod _{0\leq \lambda \leq n}~
\coprod _{I\in\I_{\lambda ,n}}~~[\bar{C}^{n-\lambda}_{r}]_{I},
\end{equation}
\setcounter{thm}{\value{equation}}

\noindent where \ $\I_{\lambda ,n}$ \ is as in \ (\ref{efive}) \ 
and \ $\bar{C}^{m}_{0}=0=\bar{C}^{0}_{r}$ \ for all \ $m,r\geq 0$.

The coface maps \ $d^{i}:E^{n-1}_{r}\to E^{n}_{r}$ \ are determined
by the cosimplicial identities and the requirement that \ 
$d^{i}|_{[\bar{C}^{n-\lambda}_{r}]_{(i_{1},\dotsc,i_{n})}}$ \ be an
isomorphism onto \ $[\bar{C}^{n-\lambda}_{r}]_{(i_{1},\dotsc,i_{n},i)}$ \ 
if \ $i>i_{n}$. \ 

The only summand in \ (\ref{enine}) \ which is not defined is thus \ 
$[\bar{C}^{n}_{r}]_{\emptyset}$, \ which we denote simply by \
$\bar{C}^{n}_{r}$. \ We require that it be an $n$-th \textit{cross-term\/}
in the sense that \ $\bar{d}_{0}|_{\bar{C}^{n}_{r}}$ \ does not 
factor through the image of any coface map \ 
$d^{i}:E^{n-1}_{r-1}\to E^{n}_{r-1}$. \ Other than that, \
$\bar{C}^{n}_{r}$ \ may be any free \Pa\ which ensures that \ 
(\ref{enine}) \ defines a $CW$-basis for a $CW$-resolution \ 
$E_{\bullet}^{n}\to\pis\U^{n}$. \ We shall call the double sequence \ 
$((\bar{C}^{n}_{r})_{n=1}^{\infty})_{r=1}^{\infty}$ \ a 
\textit{cross-term basis\/} for \ $\Edd$.
\end{defn}

Note that \ $\Ad$ \ is a retract of \ $E^{2}_{\bullet}$ \ in two
different ways (under the two coface maps \ $d^{0}$, \ $d^{1}$), \ 
corresponding to the fact that $\X$ is a retract of \ $\X\times\X$ \
in two different ways; the presence of the cross-terms \ 
$\bar{C}^{2}_{r}$ \ indicates that \ $\Ad\times\Ad$ \ is a resolution 
of \ $\pis\X^{2}$, \ but not a free one, while \ $\Ad\amalg\Ad$ \ 
is a free simplicial \Pa, \ but not a resolution.

Similarly, \ $\X\times\X$ \ embeds in \ $\X^{3}$ \ in three different 
ways, and so on\vsm .

\begin{example}\label{ect}\stepcounter{subsection}
For any \ $\Ad\to\pis\X$ \ we may set \ 
$\bar{C}^{2}_{1} = \coprod_{S^{p}_{x}\hra A^{(0)}_{0}}
           \coprod_{S^{q}_{y}\hra A^{(1)}_{0}} ~~~S^{p+q-1}_{(x,y)}$, \ 
with \ $\bar{d}_{0}|_{S^{p+q-1}_{(x,y)}}=[\iota_{x},\iota_{y}]$ \ 
(in the notation of \S \ref{dfpc}). \ The higher cross-terms \ 
$\bar{C}^{n}_{1}=0$ \ 
for \ $n\geq 3$, \ since any $k$-th order cross-term element $z$ in \
$\coprod_{j=0}^{n} A^{(j)}_{0}$ \ ($k\geq 3$) \ is a sum of elements 
of the form \ 
$z=\zeta^{\#}[\dotsc[[\iota^{r_{1}}_{(x_{1})},\iota^{r_{2}}_{(x_{2})}],
\iota^{r_{3}}_{(x_{3})}],\dotsc,\iota^{r_{k}}_{(x_{k})}]$, \ and then \ 
$$
z=d_{0}(\zeta^{\#}
[\dotsc[\iota^{r_{1}+r_{2}-1}_{(x_{1},x_{2})},s_{0}\iota^{r_{3}}_{(x_{3})}],
\dotsc,s_{0}\iota^{r_{k}}_{(x_{k})}]).
$$
\end{example}

\begin{defn}\label{dwdd}\stepcounter{subsection}
Let \ $^{h}\Wdd\to\Udd$ \ be the $\Delta$-cosimplicial augmented 
simplicial space up-to-homotopy which corresponds to \ 
$\Edd\to\pis\Udd$ \ via \S \ref{dfpa}. Thus the various (co)simplicial 
morphisms exist, and satisfy the (co)simplicial identities, only in the 
homotopy category (we may choose representatives in \ $\Ta$, \ but then
the identities are satisfied only up to homotopy). \ Each \ $\uW{n}{r}$ \ is 
homotopy equivalent to a wedge of spheres, and has a wedge summand \
$\bW{n}{r}\hra \uW{n}{r}$ \ corresponding to the $CW$-basis
free \Pa\ summand \ $\bar{E}^{n}_{r}\hra E^{n}_{r}$. \ We let \ 
$\bar{\C}^{n}_{r}$ \ denote the wedge summand of \ $\bW{n}{r}$ \
corresponding to \ $\bar{C}^{n}_{r}\hra \bar{E}^{n}_{r}$.
\end{defn}

\begin{defn}\label{drec}\stepcounter{subsection}
An simplicial space \ $\Vd\in s\Ta$ \ is called a \textit{rectification} of 
a simplicial space up-to-homotopy \ $^{h}\Wd$ \ if \ $\V_{n}\simeq \W_{n}$ \ 
for each \ $n\geq 0$, \ and the face and degeneracy maps of \ $\Vd$ \ are 
homotopic to the corresponding maps of \ $^{h}\Wd$. \ See 
\cite[\S 2.2]{DKSmH}, e.g., for a more precise definition; for our purposes 
all we require is that \ $\pis\Vd$ \ be isomorphic (as a simplicial \Pa) to \ 
$\pis(^{h}\Wd)$. \ Similarly for rectification of ($\Delta$-)cosimplicial 
objects, and so on.
\end{defn}

By considering the proof of Theorem \ref{tone}, we see that we can make the
following
%
%
\begin{assum}\label{aone}\stepcounter{subsection}
$\Edd$ \ maps monomorphically into \ $\pis\Vd(U^{\bullet}_{\Delta})$, \ 
and \ $^{h}\Wdd\to\Udd$ \ can be rectified so as to yield a strict 
$\Delta$-co\-simp\-li\-cial augmented simplicial space \ $\Wdd\to\Udd$ \ 
realizing \ $\Edd\to\pis\Udd$.
\end{assum}

\begin{defn}\label{dwud}\stepcounter{subsection}
Now assume that \ $\pis\X$ \ is an abelian \Pa\ (Def.\ \ref{dapa}) \ -- \ 
this is the necessary \Pa\ condition in order for $\X$ to be an 
$H$-space \ -- \ and let \ $\mu:\pis\X\times\pis\X\to\pis\X$ \ 
be the morphism of \Pa s defined levelwise by the group operation (see
\cite[\S 2]{BlaHO}). \ This $\mu$ is of course associative, in the sense
that \ $\mu\circ (\mu,id)=\mu\circ (id,\mu):\pis(\X^{3})\to\pis\X$, \ so
it allows one to extend the $\Delta$-cosimplicial \Pa\ \ $\Fdd\DEF\pis(\Udd)$ \
to a full cosimplicial \Pa\ \ $\Fd$, \ defined as in \S  \ref{fahs}. \ 

Since \ $E_{\bullet}^{n}\to F^{n}=\pis\U^{n}$ \ is a free resolution of \Pa s, 
the codegeneracy maps \ $s^{j}:F^{n}\to F^{n-1}$ \ induce maps of
simplicial \Pa s \ $s^{j}_{\bullet}:E_{\bullet}^{n}\to E_{\bullet}^{n-1}$, \ 
unique up to simplicial homotopy, by the universal property of
resolutions (cf.\ \cite[I, p.\ 1.14 \& II, \S 2, Prop.\ 5]{QuH}). \
Note, however, that the individual maps \ 
$s^{j}_{r}:E^{n}_{r}\to E^{n-1}_{r}$ \ are not unique, in general; in fact,
different choices may correspond to different $H$-multiplications on $\X$.
\end{defn}

These maps \ $s^{j}$ \ make \ $\Edd\to\Fdd$ \ into a full cosimplicial
augmented simplicial \Pa\ \ $\Eud\to\Fd$, \ and thus \ $^{h}\Wud\to\Udd$ \ 
into a cosimplicial augmented simplicial space up-to-homotopy \ 
(for which we may assume by \ref{aone} that all simplicial 
identities, and all the cosimplicial identities involving only the 
coface maps, hold precisely).

%
%
\begin{prop}\label{pfour}\stepcounter{subsection}
The cosimplicial simplicial space up-to-homotopy \ 
$^{h}\Wud$ \ of \S \ref{dwud} may be rectified if and only if $\X$ is homotopy
equivalent to a loop space.
\end{prop}

\begin{proof}
If $\X$ is a loop space, it has a strictly associative $H$-multiplication \ 
$m:\X\times\X\to\X$ \ which induces $\mu$ on \ $\pis(-)$ \ (cf.\
\cite[Prop.\ 9.9]{GrayH}), so \ $\Udd$ \ extends to a
cosimplicial space \ $\Ud$ \ by Fact \ref{fahs}. \ Applying the
functorial construction of \cite[\S 2]{StoV} to \ $\Ud$ \ yields a (strict)
cosimplicial augmented simplicial space \ $\Vdd\to\Ud$, \ and since we
assumed \ $\pis\uWd{n}$ \ embeds in \ $\pis\V^{n}_{\bullet}$ \ for each $n$, \
$^{h}\Wud$ \ may also be rectified\vsm .

Conversely, if \ $\Wud$ \ is a (strict) cosimplicial simplicial space 
realizing \ $\Eud$, \ then we may apply the realization functor for
simplicial spaces in each cosimplicial dimension \ $n\geq 0$ \ to obtain \  
$\|\uWd{n}\|\simeq \U^{n}=\X^{n+1}$ \ (by \S \ref{sqss}). \ 
The realization of the codegeneracy map \ 
$\|s^{0}\|:\|\uWd{1}\|\to\|\uWd{0}\|$ \ induces \
$\mu:\pis(\X^{2})\to\pis\X$, \ so it corresponds to an $H$-space 
multiplication \ $m:\X^{2}\to\X$ \ (see \cite[Prop.\ 2.7]{BlaHO}).

The fact that \ $\|\Wud\|$ \ is a (strict) cosimplicial space means that
all composite codegeneracy maps \ 
$\|s^{0}\circ s^{j_{1}}\circ\dotsb s^{j_{n-1}}\|:\|\uWd{n}\|\to\|\uWd{0}\|$ \
are equal, \ and thus all possible composite multiplications \
$\X^{n+1}\to\X$ \ (i.e., all possible bracketings in \ (\ref{eone})) \ are 
homotopic, with homotopies between the homotopies, and so on \ -- \ in
other words, the $H$-space \ $\lra{\X,m}$ \ is an
$A_{\infty}$ space (see \cite[Def.\ 11.2]{StaH}) \ -- \ so that $\X$ is 
homotopy equivalent to loop space by \cite[Theorem 11.4]{StaH}. \ Note
that we only required that the codegeneracies of \ $^{h}\Wud$ \ 
be rectified; after the fact this ensures that the full cosimplicial
simplicial space is rectifiable.
\end{proof}

In summary, the question of whether $\X$ is a loop space reduces to the
question of whether a certain diagram in the homotopy category,
corresponding to a diagram of free \Pa s, may be rectified \ -- \ 
or equivalently, may be made $\infty$-homotopy commutative.

%
%
\section{Polyhedra and higher homotopy operations}
\label{cphho}

As in \cite[\S 4]{BlaHH}, there is a sequence of higher homotopy
operations which serve as obstructions to such a rectification, and these
may be described combinatorially in terms of certain polyhedra, 
as follows:

\begin{defn}\label{dfcf}\stepcounter{subsection}
The $N$-\textit{permutohedron} \ $\Pe{N}$ \ is defined to be
the convex hull in \ $\R^{N}$ \ of the points \ 
$p_{\sigma}=(\sigma(1),\sigma(2),\dotsc,\sigma(N))$, \ where $\sigma$
ranges over all permutations \ $\sigma\in\Sigma_{N}$ \ 
(cf.\ \cite[\S 9]{ZiegL}). It is \ ($N-1$)-dimensional.

For any two integers \ $0\leq n< N$, \ the corresponding 
($N,n$)-\textit{face-codegeneracy polyhedron} \ $\PP{N}{n}$ \ is a 
quotient of the $N$-permutohedron \ $\Pe{N}$ \ obtained by identifying 
two vertices \ $p_{\sigma}$ \ and \ $p_{\sigma'}$ \ to a single vertex \ 
$\bar{p}_{\sigma}= \ \bar{p}_{\sigma'}$ \ of \ $\PP{N}{n}$ \ whenever \ 
$\sigma=(i,i+1)\sigma'$, \ where \ $(i,i+1)$ \ is an adjacent
transposition and \ $\sigma(i),\sigma(i+1)>n$. \ 

Since each facet $A$ of \ $\Pe{N}$ \ is uniquely determined by its
vertices (see below), the facets in the quotient \
$\PP{N}{n}$ \ are obtained by collapsing those of \ $\Pe{N}$ \ accordingly.
\end{defn}

Note that \ $\PP{N}{N-1}$ \ is the $N$-permutohedron \ $\Pe{N}$, \ and
in fact the quotient map \ $q:\Pe{N}\epic\PP{N}{n}$ \ is homotopic to a
homeomorphism (though not a combinatorial isomorphism, of course) for \ 
$n\geq 1$. \ On the other hand, \ $\PP{N}{0}$ \ is a single point. \ 
For non-trivial examples of face-codegeneracy polyhedra, see 
Figures \ref{fig3} \& \ref{fig4} below.

\begin{fact}\label{ffcf}\stepcounter{subsection}
From the description of the facets of the permutohedron given in 
\cite{GGupA}, we see that \ $\PP{N}{n}$ \ has an edge connecting 
a vertex \ $p_{\sigma}$ \ to any vertex of the form \ 
$p_{(i,i+1)\sigma}$ \ (unless \ $\sigma(i),\sigma(i+1)>n$, \ 
in which case the edge is degenerate).

More generally, let \ $\bar{p}_{\sigma}$ \ be any vertex of \ $\PP{N}{n}$. \ 
The facets of \ $\PP{N}{n}$ \ containing \ $\bar{p}_{\sigma}$ \ are
determined as follows:

Let \ $\bP=\lra{1,2,\dotsc,\ell_{1}~|\ \ell_{1}+1,\dotsc,\ell_{2}~|~\dotsc~
|\ \ell_{i-1}+1,\dotsc,\ell_{i}~|\ \dotsc~|\ \ell_{r-1}+1,\dotsc,N~}$ \ 
be a partition of \ $1,\dotsc,N$ \ into $r$ consecutive blocs, subject to
the condition that for each \ $1\leq j<r$ \ at least one of \ 
$\sigma(\ell_{i})$, \ $\sigma(\ell_{i+1})$ \ is \ $\leq n$. \ Denote
by \ $n_{i}$ \ the number of $j$'s  in the $i$-th bloc (i.e., \ 
$\ell_{i-1}+1\leq j\leq \ell_{i}$) \ such that \ $\sigma(j)\leq n$.\hsm 
Then \ $\PP{N}{n}$ \ will have a subpolyhedron \ $Q(\bP)$ \
(containing \ $p_{\sigma}$) \ which is isomorphic to the product \ 
$$
\PP{\ell_{1}}{n_{1}}\times\PP{\ell_{2}-\ell_{1}}{n_{2}}\times\cdots\times 
\PP{\ell_{i}-\ell_{i-1}}{n_{i}}\times\cdots\times\PP{N-\ell_{r-1}}{n_{r}}.
$$

\noindent This follows from the description of the facets of the 
$N$-permutohedron in \cite[\S 4.3]{BlaHH}.

We denote by \ $(\PP{N}{n})^{(k)}$ \ the union of all facets of \
$\PP{N}{n}$ \ of dimension $\leq k$. \ In particular, for \ $n\geq 1$ \ 
we have \ $\partial\PP{N}{n}\DEF(\PP{N}{n})^{(N-2)}=\bS{N-2}$, \ 
since the homeomorphism \ $\tilde{q}:\Pe{N}\to\PP{N}{n}$ \ preserves \ 
$\partial\Pe{N}$.
\end{fact}

\subsection{Factorizations}
\label{sfac}\stepcounter{thm}

Given a cosimplicial simplicial object \ $\Eud$ \ as in \S \ref{dwud},
any composite face-codegeneracy map \
$\psi:E^{n+k}_{m+\ell}\to E^{k}_{\ell}$ \ has a (unique) canonical
factorization of the form \ $\psi=\phi\circ\theta$, \ where \ 
$\theta:E^{n+k}_{m+\ell}\to E^{k}_{m+\ell}$ \ may be written \ 
$\theta=s^{j_{1}}\circ s^{j_{2}}\circ\dotsc s^{j_{n}}$ \ for \ 
$0\leq j_{1}<j_{2}<\dotsc<j_{n}<n+k$ \ and \ 
$\phi:E^{k}_{m+\ell}\to E^{k}_{\ell}$ \ may be written \ 
$\phi=d_{i_{1}}\circ d_{i_{2}}\circ\cdots d_{i_{n}}$ \ for \ 
$0\leq i_{1}<i_{2}<\dotsc<i_{n}\leq m+\ell$.

Let \ $\D(\psi)$ \ denote the set of all possible (not necessarily
canonical) factorizations of $\psi$ as a composite of face and
codegeneracy maps: \ $\psi=\lambda_{n+m}\circ\dotsb\circ\lambda_{1}$. \ 
We define recursively a bijective correspondence between \ $\D(\psi)$ \ 
and the vertices of an \ ($n+m$)-permutohedron \ $\Pe{n+m}$, \ as
follows (compare \cite[Lemma 4.7]{BlaHH}):

The canonical factorization \ 
$\psi=d_{i_{1}}\circ d_{i_{2}}\circ\cdots d_{i_{n}}\circ 
s^{j_{1}}\circ s^{j_{2}}\circ\cdots s^{j_{n}}$ \
corresponds to the vertex \ $p_{id}$.\hsm Next, assume that
the factorization \ 
$\psi=\lambda_{n+m}\circ\dotsc\circ\lambda_{1}$ \ corresponds to \ 
$p_{\sigma}$. \ Then the factorization corresponding to \ 
$p_{\sigma'}$, \ for \ $\sigma=(i,i+1)\sigma'$, \ 
is obtained from \ $\psi=\lambda_{1}\circ\dotsb\circ\lambda_{n+m}$ \
by switching \ $\lambda_{i}$ \ and  \ $\lambda_{i+1}$, \ using 
the identity \ $s^{j}\circ s^{i}=s^{i-1}\circ s^{j}$ \ for \ $i>j$ \
if \ $\lambda_{i}$ \ and  \ $\lambda_{i+1}$ \ are both codegeneracies,
and the identity \ $d_{i}\circ d_{j}=d_{j-1}\circ d_{i}$ \ for \ $i<j$ \ 
if they are both face maps. 

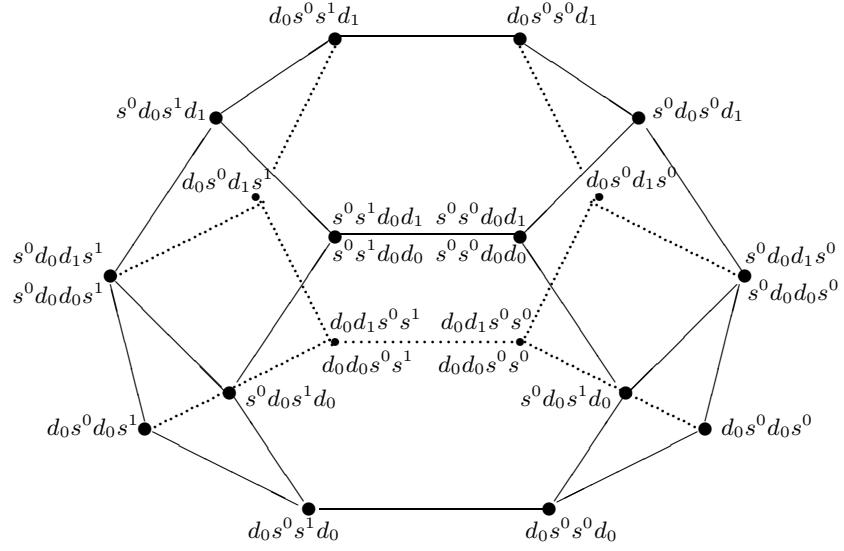
\begin{figure}[htbp]
\begin{picture}(320,220)(-50,0)
%
%
\put(80,165){\circle*{5}}
\put(42,165){\scriptsize $s^{0}d_{0}s^{1}d_{1}$}
\put(82,167){\line(3,2){43}}
\put(82,163){\line(1,-1){43}}
\put(125,195){\circle*{5}}
\put(100,201){\scriptsize $d_{0}s^{0}s^{1}d_{1}$}
\put(126,196){\line(1,0){68}}
\put(195,195){\circle*{5}}
\put(190,201){\scriptsize $d_{0}s^{0}s^{0}d_{1}$}
\put(197,193){\line(3,-2){42}}
\put(240,165){\circle*{5}}
\put(245,165){\scriptsize $s^{0}d_{0}s^{0}d_{1}$}
\put(238,163){\line(-1,-1){43}}
\put(195,120){\circle*{5}}
\put(163,125){\scriptsize $s^{0}s^{0}d_{0}d_{1}$}
\put(163,111){\scriptsize $s^{0}s^{0}d_{0}d_{0}$}
\put(125,120){\circle*{5}}
\put(124,125){\scriptsize $s^{0}s^{1}d_{0}d_{1}$}
\put(124,111){\scriptsize $s^{0}s^{1}d_{0}d_{0}$}
\put(126,121){\line(1,0){70}}
%
%
\put(40,105){\circle*{5}}
\put(3,110){\scriptsize $s^{0}d_{0}d_{1}s^{1}$}
\put(3,95){\scriptsize $s^{0}d_{0}d_{0}s^{1}$}
\put(42,107){\line(2,3){36}}
\put(42,103){\line(1,-1){42}}
\put(85,61){\circle*{5}}
\put(91,56){\scriptsize $s^{0}d_{0}s^{1}d_{0}$}
\put(88,65){\line(2,3){35}}
%
%
\put(280,105){\circle*{5}}
\put(280,110){\scriptsize $s^{0}d_{0}d_{1}s^{0}$}
\put(281,96){\scriptsize $s^{0}d_{0}d_{0}s^{0}$}
\put(279,107){\line(-2,3){36}}
\put(277,102){\line(-1,-1){42}}
\put(235,61){\circle*{5}}
\put(195,56){\scriptsize $s^{0}d_{0}s^{1}d_{0}$}
\put(231,65){\line(-2,3){35}}
%
%
\put(115,17){\circle*{5}}
\put(92,6){\scriptsize $d_{0}s^{0}s^{1}d_{0}$}
\put(119,17){\line(1,0){86}}
\put(113,20){\line(-2,3){26}}
\put(206,17){\circle*{5}}
\put(197,6){\scriptsize $d_{0}s^{0}s^{0}d_{0}$}
\put(208,20){\line(2,3){26}}
%
%
\put(53,47){\circle*{5}}
\put(16,45){\scriptsize $d_{0}s^{0}d_{0}s^{1}$}
\put(56,46){\line(2,-1){55}}
\put(53,50){\line(-1,4){13}}
%
%
\put(265,47){\circle*{5}}
\put(271,45){\scriptsize $d_{0}s^{0}d_{0}s^{0}$}
\put(262,46){\line(-2,-1){53}}
\put(265,50){\line(1,4){13}}
%
%
\put(95,135){\circle*{3}}
\put(67,138){\scriptsize $d_{0}s^{0}d_{1}s^{1}$}
\multiput(97,133)(1,-2){28}{\circle*{.5}}
\multiput(102,146)(1,2){24}{\circle*{.5}}
\put(125,80){\circle*{3}}
\put(123,85){\scriptsize $d_{0}d_{1}s^{0}s^{1}$}
\put(120,69){\scriptsize $d_{0}d_{0}s^{0}s^{1}$}
\multiput(126,80)(3,0){23}{\circle*{.5}}
\put(195,80){\circle*{3}}
\put(165,85){\scriptsize $d_{0}d_{1}s^{0}s^{0}$}
\put(164,69){\scriptsize $d_{0}d_{0}s^{0}s^{0}$}
\multiput(197,82)(1,2){28}{\circle*{.5}}
\put(225,135){\circle*{3}}
\put(220,139){\scriptsize $d_{0}s^{0}d_{1}s^{0}$}
\multiput(218,146)(-1,2){24}{\circle*{.5}}
%
%
\multiput(42,105)(2,1){29}{\circle*{.5}}
\multiput(57,48)(2,1){34}{\circle*{.5}}
%
%
\multiput(278,105)(-2,1){29}{\circle*{.5}}
\multiput(261,48)(-2,1){34}{\circle*{.5}}
\end{picture}
\caption[fig3]{The face-codegeneracy polyhedron 
$\protect\PP{4}{2}(d_{0}d_{1}s^{0}s^{1})$}
\label{fig3}
\end{figure}

Passing to the quotient face-codegeneracy polyhedron,
we see that the vertices of \ $\PP{n+m}{n}$ \ are now identified with
factorizations of $\psi$ of the form
\setcounter{equation}{\value{thm}}\stepcounter{subsection} 
\begin{equation}\label{eeleven}
E^{n+k}_{m+\ell} \xra{s^{j_{n_{t}}^{t}}} E^{n+k-1}_{m+\ell}\dotsc 
E^{n_{t}+1}_{m+\ell} \xra{s^{j_{1}^{t}}} E^{n_{t}}_{m+\ell}
\xra{\theta_{t}} E^{n_{t}}_{m_{t}}\dotsc E^{n_{1}}_{m_{1}} 
\xra{s^{j_{n_{1}}^{0}}}\dotsc E^{n+1}_{m_{1}} \xra{s^{j_{n_{0}}^{0}}} 
E^{n}_{m_{1}} \xra{\theta_{0}} E^{n}_{m}
\end{equation}
\setcounter{thm}{\value{equation}}
\noindent where \ $\theta_{i}$ \ is a composite of face maps (i.e., we
do not distinguish the different ways of decomposing \ $\theta_{i}$ \
as \ $d_{k_{1}}\circ\dotsc d_{k_{r}}$). \ 
The collection of such factorizations of $\psi$ will be denoted by \
$D(\psi)\msim$, \ where \ $\sim$ \ is the obvious equivalence relation
on \ $D(\psi)$. \ We shall denote the
face-codegeneracy polyhedron \ $\PP{n+m}{n}$ \ with its vertices so
labelled by \ $\PP{n+m}{n}(\psi)$. \ An example for \
$\psi=d_{0}d_{1}s^{0}s^{1}$ \ appears in Figure \ref{fig3}.

%
%

\subsection{Notation}
\label{ncc}\stepcounter{thm}
For \ $\psi:E^{n+k}_{m+\ell}\to E^{k}_{\ell}$ \ as above, we denote by \ 
$\C(\psi)$ \ the collection of all composite face-codegeneracy maps \ 
$\rho:E^{n(\rho)+k(\rho)}_{m(\rho)+\ell(\rho)}\to E^{k(\rho)}_{\ell(\rho)}$ \ 
such that $\rho$ is of the form \ $\rho=\xi_{t}\circ\dotsb\circ\xi_{s}$ \ 
($1\leq s\leq t\leq\nu$) \ for some decomposition \ 
$\psi=\xi_{\nu}\circ\dotsb\circ\xi_{1}=
\theta_{0}\circ s^{j_{n_{0}}^{0}}\circ\dotsb\circ s^{j_{n_{1}}^{0}}
\circ\theta_{1}\circ\dotsb\circ\theta_{t}\circ 
s^{j_{1}^{t}}\circ\dotsb \circ s^{j_{n_{t}}^{t}}$ \ of \ (\ref{eeleven}). \ 
That is, we allow only those subsequences \
$\lambda_{b},\dotsc,\lambda_{a}$ \ of \ a factorization \ 
$\psi=\lambda_{n+m}\circ\cdots\circ\lambda_{1}$ \ in \ $\D(\psi)$ \
which are compatible with the equivalence relation $\sim$ in the sense
that \ $\lambda_{b+1}$ \ and \ $\lambda_{b}$ \ are not both face maps,
and similarly for \ $\lambda_{a-1}$ \ and \ $\lambda_{a}$.\hsm
Such a $\rho$ will be called \textit{allowable}.

\subsection{Higher homotopy operations}
\label{shho}\stepcounter{thm}

Given a cosimplicial simplicial space up-to-homotopy \ $^{h}\Wud$ \ as
in \S \ref{scsc}, we now define a certain sequence of higher homotopy 
operations. \ First recall that the \textit{half-smash} of two spaces \
$\X,\Y\in\Ta$ \ is  \ $\X\ltimes\Y\DEF (\X\times\Y)/(\X\times \{\ast\})$; \ 
if \ $\X$ is a suspension, there is a (non-canonical) homotopy 
equivalence \ $\X\ltimes\Y \simeq \X\wedge \Y\vee \X$.

\begin{defn}\label{dcc}\stepcounter{subsection}
Given a composite face-codegeneracy map \ 
$\psi:\uW{n+k}{m+\ell}\to \uW{k}{\ell}$ \ as above, a 
\textit{compatible collection for} $\C(\psi)$ \ and \ $^{h}\Wud$ \ is a set \
$\{g^{\rho}\}_{\rho\in\C(\psi)}$ \ of maps \
$g^{\rho}:\PP{n(\rho)+m(\rho)}{m(\rho)}(\rho)\ltimes
\uW{n(\rho)+k(\rho)}{m(\rho)+\ell(\rho)}\to\uW{k(\rho)}{\ell(\rho)}$ \ 
for each \ $\rho\in\C(\psi)$, \ satisfying the following condition\vsm :

Assume that for such a \ $\rho\in\C(\psi)$ \ we have
some decomposition \ 
$$
\rho=\xi_{\nu}\circ\dotsb\circ\xi_{1}=
\theta_{0}\circ s^{j_{n_{0}}^{0}}\circ\dotsc\circ s^{j_{n_{1}}^{0}}
\circ\theta_{1}\circ\dotsb\circ\theta_{t}\circ 
s^{j_{1}^{t}}\circ\dotsb \circ s^{j_{n_{t}}^{t}}
$$
\noindent in \ $\D(\rho)\msim$, \ as in \ (\ref{eeleven}), \ and let \ 
$$
\bP=\lra{1,\dotsc,\ell_{1}~|\ \dotsc~|\ 
\ell_{i-1}+1,\dotsc,\ell_{i}~|\ \dotsc~|\ \ell_{r-1}+1,\dotsc,\nu~}
$$
be a partition of \ $(1,\dotsc,\nu)$ \ as in \S \ref{ffcf}, \ 
yielding a sequence of composite face-codegeneracy maps \ 
$\rho_{i}\in\C(\rho)\subseteq\C(\psi)$ \ for \ $i=1,\dotsc,r$.

Let \ $Q(\bP)\cong \PP{\ell_{1}}{n_{1}}(\rho_{1})\times\cdots\times
\PP{\ell_{i}-\ell_{i-1}}{n_{i}}(\rho_{i})\times\dotsb\times
\PP{\nu-\ell_{r-1}}{n_{r}}(\rho_{r})$ \ be the corresponding sub-polyhedron
of \ $\PP{n(\rho)+m(\rho)}{m(\rho)}(\rho)$. \ Then we require that \ 
$g^{\rho}|_{Q(\bP)\ltimes\uW{n(\rho)+k(\rho)}{m(\rho)+\ell(\rho)}}$ \
be the composite of the corresponding maps \ $g^{\rho_{i}}$ \ in the
sense that
\setcounter{equation}{\value{thm}}\stepcounter{subsection} 
\begin{equation}\label{etwelve}
g^{\rho}(x_{1},\dotsc,x_{r},w)=g^{\rho_{1}}(x_{1},g^{\rho_{2}}(x_{2},
\dotsc,g^{\rho_{r}}(x_{r},w)\dotsc))
\end{equation}
\setcounter{thm}{\value{equation}}
for \ $x_{i}\in \PP{\ell_{i}-\ell_{i-1}}{n_{i}}(\rho_{i})$ \ and \ 
$w\in\uW{n(\rho)+k(\rho)}{m(\rho)+\ell(\rho)}$. 

We further require that if \ $\rho=\lambda_{1}$ \ is of length $1$, then \
$g^{\rho}$ \ must be in the prescribed homotopy class of the face or
codegeneracy map \ $\lambda_{1}$. \ Thus in particular, for each 
vertex \ $\bar{p}_{\sigma}$ \ of \ $\PP{n+m}{n}(\psi)$, \ indexed 
by a factorization \ $\psi=\xi_{\nu}\circ\dotsb\circ\xi_{1}$ \ in \ 
$\D(\psi)\msim$, \ the map \ 
$g^{\rho}|_{\{\bar{p}_{\sigma}\}\times \uW{n+\ell}{m+k}}$ \ 
represents the class \ $[\xi_{\nu}\circ\dotsb\circ\xi_{1}]$.
\end{defn}

\begin{fact}\label{fim}\stepcounter{subsection}
Any compatible collection of maps \ $\{g^{\rho}\}_{\rho\in\C(\psi)}$ \
for \ $C(\psi)$ \ induces a map \ 
$f=f^{\psi}:\partial\PP{n+m}{n}\ltimes\uW{n+k}{m+\ell}\to\uW{k}{\ell}$ \ 
(since all the facets of \ $\partial\PP{n+m}{n}$ \ are products of \ 
face-codegeneracy polyhedra of the form \ 
$\PP{n(\rho)+m(\rho)}{n(\rho)}(\rho)$ \ for \ $\rho\in\C(\psi)$, \ 
and condition \ (\ref{etwelve}) \ guarantees that the maps \ 
$g^{\rho}$ \ agree on intersections).
\end{fact}

\begin{defn}\label{dhho}\stepcounter{subsection}
Given \ $^{h}\Wud$ \ as in \S \ref{dwud}, for each \ $k\geq 2$ \ and each 
composite face-codegeneracy map \ $\psi:\uW{n+k}{m+\ell}\to \uW{k}{\ell}$, \ 
the $k$-{\em th order homotopy operation\/} 
associated to \ $^{h}\Wud$ \ and $\psi$ \ is 
a subset \ $\lra{\psi}$ \ of the track group \ 
$[\Sigma^{n+m-2}\uW{n+k}{m+\ell},\uW{k}{\ell}]$, \ defined as follows:

Let \ 
$S\subseteq[\partial\PP{n+m}{n}\ltimes\uW{n+k}{m+\ell},\ \uW{k}{\ell}]$ \ 
be the set of homotopy classes of maps \ 
$f=f^{\psi}:\partial\PP{n+m}{n}\ltimes\uW{n+k}{m+\ell}\to\uW{k}{\ell}$ \ 
which are induced as above by some compatible collection \ 
$\{ g^{\rho} \} _{\rho\in\C(\psi)}$ \ for \ $\C(\psi)$.

Now choose a splitting 
\setcounter{equation}{\value{thm}}\stepcounter{subsection} 
\begin{equation}\label{ethirteen}
\partial\PP{n+m}{n}(\psi)\ltimes\uW{n+k}{m+\ell}\cong
\bS{n+m-2}\ltimes\uW{n+k}{m+\ell}\simeq
(\bS{n+m-2}\wedge\uW{k}{\ell})\vee\uW{k}{\ell}
\end{equation}
\setcounter{thm}{\value{equation}}
and let \ 
$\lra{\psi}\subseteq[\Sigma^{n+m-2}\uW{n+k}{m+\ell},\uW{k}{\ell}]$ \ 
be the image of the subset $S$ under the resulting projection.
\end{defn}

It is clearly a necessary condition in order for the subset \ 
$\lra{\psi}$ \ to be non-empty that all the lower order operations \
$\lra{\rho}$ \  {\em vanish\/} \ (i.e., contain the null class)
for all \ $\rho\in\C(\psi)\setminus\{\psi\}$ \ -- \ because otherwise 
the various maps \ 
$g^{\rho}:\PP{n(\rho)+m(\rho)}{m(\rho)}(\rho)\ltimes
\uW{n(\rho)+k(\rho)}{m(\rho)+\ell(\rho)}\to\uW{k(\rho)}{\ell(\rho)}$ \
cannot even extend over the interior of \
$\PP{n(\rho)+m(\rho)}{m(\rho)}(\rho)$.\hsm 
A {\em sufficent\/} condition is that the operations \ $\lra{\rho}$ \ 
vanish {\em coherently\/}, in the sense that the choices of compatible
collections for the various $\rho$ be consistent on common subpolyhedra 
(see \cite[\S 5.7]{BlaHH} for the precise definition, and 
\cite[\S 5.9]{BlaHH} for the obstructions to coherence).

On the other hand, if \ $^{h}\Wud$ \ is the cosimplicial simplicial 
space up-to-homotopy of \S \ref{dedd} (corresponding to the
cosimplicial simplicial \Pa\ \ $\Edd$ \ with the $CW$-basis \ 
$\{\bar{E}_{r}^{n}\}_{r,n=0}^{\infty}$), \ 
then the vanishing of the homotopy operation \ 
$\lra{\psi|_{\bar{\C}^{n}_{r}}}$ \ -- \ with $\psi$ 
restricted to the \ ($n,r$)-cross-term \ -- \ implies the vanishing 
of \ $\lra{\psi}$, \ for \ any $\psi:\uW{n+k}{m+\ell}\to \uW{k}{\ell}$ \ 
(assuming lower order vanishing). \ This is because outside of the 
wedge summand \ $\bar{\C}^{n}_{r}$, \ the map $\psi$ is 
determined by the maps \ $\rho\in\C(\psi)$ \ and the coface and 
degeneracy maps of \ $^{h}\Wud$, \ which we may assume to 
$\infty$-homotopy commute by induction and \ref{aone} 
respectively. 

We may thus sum up the results of this section, combined with Proposition 
\ref{pfour}, in:

%
%
\begin{thm}\label{tthree}\stepcounter{subsection}
A space \ $\X\in\Ta$, \ for which \ $\pis\X$ \ is an abelian \Pa, is homotopy 
equivalent to a loop space if and only if all the higher homotopy operations \ 
$\lra{\psi|_{\bar{\C}^{n}_{r}}}$ \ defined above vanish coherently.
\end{thm}

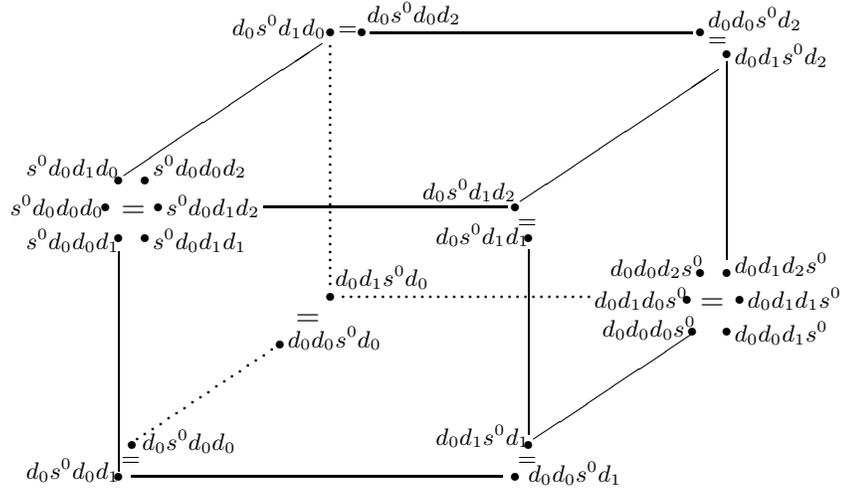
\begin{figure}[htbp]
\begin{picture}(300,185)(-40,0)
%
%
\put(46,106){$=$}
\put(4,107){\scriptsize $s^{0}d_{0}d_{0}d_{0}$}
\put(40,110){\circle*{3}}
\put(10,94){\scriptsize $s^{0}d_{0}d_{0}d_{1}$}
\put(45,98){\circle*{3}}
\put(58,94){\scriptsize $s^{0}d_{0}d_{1}d_{1}$}
\put(55,98){\circle*{3}}
\put(63,107){\scriptsize $s^{0}d_{0}d_{1}d_{2}$}
\put(60,110){\circle*{3}}
\put(58,122){\scriptsize $s^{0}d_{0}d_{0}d_{2}$}
\put(55,120){\circle*{3}}
\put(10,122){\scriptsize $s^{0}d_{0}d_{1}d_{0}$}
\put(45,120){\circle*{3}}
\put(45,93){\line(0,-1){83}}
\put(100,110){\line(1,0){92}}
%
%
\put(46,12){\scriptsize $=$}
\put(54,18){\scriptsize $d_{0}s^{0}d_{0}d_{0}$}
\put(50,20){\circle*{3}}
\put(10,7){\scriptsize $d_{0}s^{0}d_{0}d_{1}$}
\put(45,8){\circle*{3}}
\put(50,8){\line(1,0){140}}
\put(196,102){\scriptsize $=$}
\put(160,113){\scriptsize $d_{0}s^{0}d_{1}d_{2}$}
\put(195,110){\circle*{3}}
\put(165,96){\scriptsize $d_{0}s^{0}d_{1}d_{1}$}
\put(200,98){\circle*{3}}
\put(200,94){\line(0,-1){70}}
\put(196,12){\scriptsize $=$}
\put(165,20){\scriptsize $d_{0}d_{1}s^{0}d_{1}$}
\put(200,20){\circle*{3}}
\put(200,6){\scriptsize $d_{0}d_{0}s^{0}d_{1}$}
\put(195,8){\circle*{3}}
%
%
\put(47,122){\line(3,2){75}}
\put(128,175){\scriptsize $=$}
\put(88,175){\scriptsize $d_{0}s^{0}d_{1}d_{0}$}
\put(125,176){\circle*{3}}
\put(139,180){\scriptsize $d_{0}s^{0}d_{0}d_{2}$}
\put(137,176){\circle*{3}}
\put(140,176){\line(1,0){122}}
\put(268,171){\scriptsize $=$}
\put(268,178){\scriptsize $d_{0}d_{0}s^{0}d_{2}$}
\put(265,176){\circle*{3}}
\put(278,164){\scriptsize $d_{0}d_{1}s^{0}d_{2}$}
\put(275,168){\circle*{3}}
\put(197,112){\line(3,2){75}}
\put(275,165){\line(0,-1){75}}
%
%
\put(202,22){\line(3,2){60}}
%
%
\put(265,71){$=$}
\put(225,72){\scriptsize $d_{0}d_{1}d_{0}s^{0}$}
\put(260,75){\circle*{3}}
\put(228,61){\scriptsize $d_{0}d_{0}d_{0}s^{0}$}
\put(262,63){\circle*{3}}
\put(278,59){\scriptsize $d_{0}d_{0}d_{1}s^{0}$}
\put(275,63){\circle*{3}}
\put(283,72){\scriptsize $d_{0}d_{1}d_{1}s^{0}$}
\put(280,75){\circle*{3}}
\put(278,85){\scriptsize $d_{0}d_{1}d_{2}s^{0}$}
\put(275,85){\circle*{3}}
\put(231,84){\scriptsize $d_{0}d_{0}d_{2}s^{0}$}
\put(265,85){\circle*{3}}
%
%
\multiput(125,171)(0,-3){31}{\circle*{.5}}
\put(112,65){$=$}
\put(127,80){\scriptsize $d_{0}d_{1}s^{0}d_{0}$}
\put(125,76){\circle*{3}}
\multiput(52,22)(3,2){18}{\circle*{.5}}
\put(106,58){\circle*{3}}
\put(109,57){\scriptsize $d_{0}d_{0}s^{0}d_{0}$}
\multiput(130,76)(3,0){32}{\circle*{.5}}
\end{picture}
\caption[fig4]{The face-codegeneracy polyhedron 
$\protect\PP{4}{1}(d_{0}d_{1}d_{2}s^{0})$}
\label{fig4}
\end{figure}

\begin{remark}\label{rhm}\stepcounter{subsection}
As observed in \S \ref{scsc}, \ for any \ $\X\in\Ta$ \ the space \ $J\X$ \ is
the colimit of the $\Delta$-cosimplicial space \ 
$\U(X)^{\bullet}_{\Delta}$, \ and in fact the $n$-th stage of the 
James construction, \ $J_{n}\X$, \ is the (homotopy) colimit of the \ 
$(n-1)$-coskeleton of \ $\Udd$. \ Thus if we think of the sequence of 
higher homotopy operations ``in the simplicial direction'' as obstructions
to the validity of the identity \cite[Thm.\ 5.7($\ast$)]{BlaL} 
(up to $\infty$-homotopy commutativity), \ then the $n$-th cosimplicial 
dimension corresponds to verifying this identity for \ 
$f\circ i_{A}:\A\to F\B$ \ of James filtration \ $n+1$ \ (cf.\
\cite[\S 2]{JamFH})\vsm .

In particular, if we fix \ $k=\ell=0$, \ $n=1$ \ and proceed by
induction on $m$, we are computing the obstructions for the
existence of an $H$-multiplication on $\X$, \ as in
\cite{BlaHO}. (Thus if $\X$ is endowed with an $H$-space structure to
begin with, they must all vanish.) \ Observe that the
face-codegeneracy polyhedron \ $\PP{n}{1}$ \ is an \ ($n-1$)-cube, as 
in Figure \ref{fig4}, rather than the ($n-1$)-simplex we had in 
\cite[\S 4]{BlaHO} \ -- \ so the homotopy operations we obtain here 
are more complicated. This is because they take value in the homotopy 
groups of spheres, rather than those of the space $\X$.
\end{remark}

%
%

As a corollary to Theorem \ref{tthree} we may deduce the following
result of Hilton  (cf.\ \cite[Theorem C]{HilLS}):
%
%
\begin{cor}\label{cfour}\stepcounter{subsection}
If \ $\lra{\X,m}$ is a \ ($p-1$)-connected $H$-space with \ 
$\pi_{i}\X=0$ \ for \ $i\geq 3p$, \ then $\X$ is a loop space, up to
homotopy.
\end{cor}

\begin{proof}
Choose a $CW$-resolution of \ $\pis\X$ \ which is
($p-1$)-connected in each simplicial dimension, and let \ $\Eud$ \ be
as in \S \ref{dedd}. By definition 
of the cross-term \Pa s \ $C^{n}_{r}$ \ in \S \ref{dedd}, \ they must
involve Whitehead products of elements from {\em all\/} lower order 
cross-terms; but since $\X$ is an $H$-space by assumption, 
all obstructions of the form \ $\lra{\psi|_{\bar{\C}^{1}_{r}}}$ \ 
vanish  (see \S \ref{rhm}). \ Thus, the lowest dimensional obstruction
possible is a third-order operation \ 
$\lra{\psi|_{\bar{\C}^{2}_{r}}}$ \ ($r\geq 2$), \ 
which involves a triple Whitehead product and thus takes value in \ 
$\pi_{i}\uW{k}{\ell}$ \ for \ $i\geq 3p$. \ If we apply the \
($3p-1$)-Postnikov approximation functor to \ $^{h}\Wud$ \ in each
dimension, to obtain \ $^{h}\Z_{\bullet}^{\bullet}$, \ all 
obstructions to rectification vanish, and from the spectral sequence 
of \S \ref{sqss} we see that the obvious map \ 
$\X=\|\uW{1}{\bullet}\|\to\|\Z^{1}_{\bullet}\|$ \
induces an isomorphism in \ $\pi_{i}$ \ for \ $i<3p$. \ Since \ 
$\|\Z^{1}_{\bullet}\|$ \ is a loop space by Theorem \ref{tthree}, so is its
($3p-1$)-Postnikov approximation, namely $\X$.
\end{proof}

\begin{example}\label{ess}\stepcounter{subsection}
The $7$-sphere is an $H$-space (under the Cayley multiplication, 
for example), but none of the $120$ possible $H$-multiplications on \ 
$\bS{7}$ \ are homotopy-associative; the first obstruction to
homotopy-associativity is a certain ``separation element'' in \
$\pi_{21}\bS{7}$ \ (cf.\ \cite[Theorem 1.4 and Corollary 2.5]{JamM2}).

Since \ $\pis\bS{7}$ \ is a free \Pa, it has a very simple
$CW$-resolution \ $\Ad\to\pis\bS{7}$, \ with \
$\bar{A}_{0}\cong\pis\bS{7}$ \ (generated by \ $\iota^{7}$), 
and \ $\bar{A}_{r}=0$ \ for \ $r\geq 1$.\hsm
A cross-term basis (\S \ref{dedd}) for the cosimplicial simplicial 
\Pa\ \ $\Eud$ \ of \S \ref{dwud} \ is then given in dimensions $<24$ \ by\vsm :

\noindent $\bullet$\hsm $\bar{C}^{1}_{1} \cong\pis\bS{13}$, \ 
with \ $\bar{d}_{0}\iota^{13}=[d^{0}\iota^{7}, d^{1}\iota^{7}]$\vsm ;

\noindent $\bullet$\hsm $\bar{C}^{2}_{2}\cong\pis\bS{19}$, \ with \ 
$\bar{d}_{0}\iota^{19}=
[d^{0}\iota^{13}, s_{0}d^{2}d^{1}\iota^{7}] -
[d^{1}\iota^{13}, s_{0}d^{2}d^{0}\iota^{7}] +
[d^{2}\iota^{13}, s_{0}d^{1}d^{0}\iota^{7}]$\vsm ;

\noindent $\bullet$\hsm $\bar{C}^{n}_{r}$ \ is at least $24$-connected
for all other \ $n$, $r$\vsm .

We set \ $s^{j}_{r}|_{\bar{C}^{n}_{r}}=0$ \ for all \ $n\leq 2$; \
this determines \ $\Eud$ \ in degrees $\leq 21$ \ and cosimplicial
dimensions $\leq 2$\vs .

By Remark \ref{rhm}, the two secondary operations \ 
$\lra{d_{0}s^{0}|_{\bar{\C}^{1}_{1}}}$ \ and \ 
$\lra{d_{1}s^{0}|_{\bar{\C}^{1}_{1}}}$ \ must vanish; \ 
on the other hand, by Corollary \ref{cfour} all obstructions to \ 
$\bS{7}$ \ being a loop space are in degrees $\geq 21$, \ so 
the only relevant cross-term is \ $\bar{C}^{2}_{2}$, \ with three
possible third-order operations \ $\lra{\psi|_{\bar{\C}^{2}_{2}}}$, \  
for \ $\psi=d_{0}d_{1}s^{0}s^{1}$, \ 
$d_{0}d_{2}s^{0}s^{1}$, \ or \ $d_{1}d_{2}s^{0}s^{1}$. \ The
corresponding face-codegeneracy polyhedra \ $P^{4}_{2}(\psi)$ \ is as
in Figure \ref{fig4}.

It is straightforward to verify that the operations \ 
$\lra{\psi|_{\bar{\C}^{2}_{2}}}$ \ are trivial for \ 
$\psi=d_{0}d_{2}s^{0}s^{1}$ \ or \ $d_{1}d_{2}s^{0}s^{1}$ \ (in fact, many
of the maps \ $g^{\rho}$, \ for \ $\rho\in C(\psi)$, \ may be chosen
to be null). \ On may also show that there is a compatible collection \ 
$\{ g^{\rho}\}_{\rho\in C(\varphi)}$ \ for \ $\varphi=d_{0}d_{1}s^{0}s^{1}$, \ 
in the sense of \S \ref{dcc}, so that the corresponding subset \ 
$\lra{\varphi|_{\bar{\C}^{2}_{2}}}\subseteq \pi_{21}\bS{7}$ \ is
non-empty; in fact, it contains the only possible obstruction to the
$21$-Postnikov approximation for \ $\bS{7}$ \ to be a loop space. \

The existence of the tertiary operation \ 
$\lra{\varphi|_{\bar{\C}^{2}_{2}}}$ \ corresponds to the fact that 
the element \ 
$[[\iota^{7},\iota^{7}],\iota^{7}]-[[\iota^{7},\iota^{7}],\iota^{7}]+
[[\iota^{7},\iota^{7}],\iota^{7}]\in \pi_{21}\bS{7}$ \ is trivial ``for
three different reasons'': \ because of the Jacobi identity, because all
Whitehead products vanish in \ $\pis\bS{7}$, \ and because of the 
linearity of the Whitehead product \ -- \ i.e., \ $[0,\alpha]=0$\vsm .

On the other hand, we know that there {\em is\/} a $3$-primary obstruction to 
the homotopy-associativity of any $H$-multiplication on \ $\bS{7}$, \ 
namely the element \ $\sigma_{14}^{\#}\tau_{7}\in\pi_{21}\bS{7}$ \ 
(see \cite[Theorem 2.6]{JamM2}). \ 
We deduce that \ 
$0\not\in \lra{\varphi|_{\bar{\C}^{2}_{2}}}$, \ and in fact (modulo $3$)
this tertiary operation consists exactly of the elements \ 
$\pm\sigma_{14}^{\#}\tau_{7}$.

For a detailed calculation of such higher order operations using
simplicial resolutions of \Pa s, see \cite[\S 4.13]{BlaHO}.
\end{example}

\begin{remark}\label{rsa}\stepcounter{subsection}
Our approach to the question of whether $\X$ is a loop space is
clearly based on, and closely related to, the classical approaches of
Sugawara and Stasheff (cf.\ \cite{StaH1,StaH2,SugG}. \ One might
wonder why Stasheff's associahedra \ $K_{i}$ \ (cf.\
\cite[\S 2,6]{StaH1}) do not show up among the face-codegeneracy
polyhedra we describe above.  Apparently this is because we do not
work directly with the space $\X$, but rather with its simplicial
resolution, which may be thought of as a ``decomposition'' of $\X$
into wedges of spheres.
\end{remark}

\end{document}